\documentclass{article}

\usepackage[english]{babel}

\usepackage[a4paper,top=2cm,bottom=2cm,left=3cm,right=3cm,marginparwidth=1.75cm]{geometry}

\usepackage{graphicx,amssymb,lineno}
\usepackage{amsmath}
 \allowdisplaybreaks
\usepackage[colorlinks=true, allcolors=blue]{hyperref}

\begin{document}

\title{A new flow model for a thin viscous fluid film between two moving surfaces}
\author{Jos\'e M. Rodr\'{\i}guez${}^{\ref{dir1}, \ref{dir2}}$, Raquel Taboada-Vázquez${}^{\ref{dir1}, \ref{dir3}}$}
\date{}

\maketitle

{\footnotesize 
\begin{enumerate}
\item CITMAga, Santiago de Compostela, Spain. \label{dir1}
\item Department of Mathematics, Higher Technical University College of Architecture, Universidade da Coruña, A Coruña, Spain. E-mail: jose.rodriguez.seijo@udc.es. \label{dir2}
\item Department of Mathematics, School of Civil Engineering, Universidade da Coruña, A Coruña, Spain. E-mail: raquel.taboada@udc.es. \label{dir3}
\end{enumerate}
}

\begin{abstract}
We propose a two-dimensional flow model of a viscous fluid between two close moving surfaces. We show that its asymptotic behavior, when the distance between the two surfaces tends to zero, is the same as that of the the Navier-Stokes equations. 

The solutions of the new  model and Navier-Stokes equations converge to the same limit problem, that depends on the boundary conditions. If slip velocity boundary conditions are imposed on the upper and lower bound surfaces, the limit is solution of a lubrication model, but if the tractions and friction forces are known on both bound surfaces, the limit is solution of a shallow water model.

The model proposed has been obtained to be
 a valuable tool for computing viscous fluid flow between two close moving surfaces, without the need to decide a priori whether the flow is typical of a lubrication problem or it is of shallow water type, and without the enormous computational effort that would be required to solve the Navier-Stokes equations in such a thin domain.

\end{abstract}


\section{Introduction}\label{sec1}

In our previous work \cite{RodTabJMAA2021}, we have used the asymptotic expansions technique to study the behavior of a viscous fluid that flows between two very close moving surfaces. Asymptotic analysis is a mathematical tool that has been used successfully (since the pioneering works of Dean \cite{Dean1}-\cite{Dean2},  Friedrichs and Dressler \cite{FriedrichsDressler} and Goldenveizer \cite{Goldenveizer}) to obtain and justify mathematical models in solid mechanics \cite{Rigolot1972}-\cite{TutekAganovicNedelec} and fluid mechanics \cite{Cimatti1983}-\cite{PanasenkoPileckas2}, when at least one of the dimensions of the domain is much smaller than the others. Using the same mathematical technique, the authors have also proposed several new shallow water models \cite{RTV1}-\cite{RTV6} and  curved-pipe flow models \cite{CR1}-\cite{CR2}.

We have observed, in our prior article \cite{RodTabJMAA2021},  that the viscous fluid that moves between two close surfaces has two very different behaviors, depending on the boundary conditions of the problem. If the pressure differences are large in the open part of the domain boundary (that is, the region of the domain boundary between the two surfaces), then the fluid obeys equation \eqref{Reynolds_gen}, which resembles a lubrication problem. If the pressure differences are small in the mentioned region of the domain boundary, then the fluid obeys equation \eqref{ec_Vi0-v2}, which is similar to a shallow water problem.

This behavior reminds us of that observed in the works of Ciarlet et al. \cite{CiarletLodsI}-\cite{CiarletLodsIV}, where it is shown that the solution of the linearized elasticity equations in a shell converges, when the shell thickness tends to zero, to different shell models, depending on the geometry of the shell and its boundary conditions. In particular, in the work of Ciarlet and Lods \cite{CiarletLodsIII}, the authors show that the Koiter’s shell model has the same asymptotic behavior, that is, its solutions converge, when the thickness of the shell tends to zero, to the same limit problems as the linearized elasticity equations do.

In this article we intend to justify a new two-dimensional flow model of a viscous fluid between two very close moving surfaces in a similar way to what was done in the above mentioned works \cite{CiarletLodsI}-\cite{CiarletLodsIII}. We will propose the new two-dimensional model in section \ref{NewModel} and then we will justify that its asymptotic behavior is the same as that obtained in our preceding work \cite{RodTabJMAA2021} for the Navier-Stokes equations.

In the first place, we will summarize in section \ref{PreviousModels} the results presented previously in our article \cite{RodTabJMAA2021}, which will allow us to make some assumptions about the behavior of the solutions of the Navier-Stokes equations, that will be used in sections \ref{NewAssumptions}, \ref{BoundaryConditions} and \ref{NewModel} to derive a new two-dimensional flow model of a viscous fluid between two close moving surfaces. Its asymptotic behavior when the distance between the two surfaces tends to zero is justified too.

\section{Summary of the main previous results} \label{PreviousModels}

In our prior work \cite{RodTabJMAA2021} we studied the behavior of the Navier-Stokes equations in a domain bounded by two very close moving surfaces, when the distance between them tends to zero. We observed that the asymptotic behavior of the solutions of the Navier-Stokes equations in this case strongly depends on the boundary conditions. In fact, two different limit models are obtained (one similar to a lubrication model and the other similar to a shallow water model), depending on the boundary conditions in the original problem.

The two models presented in the preceding article \cite{RodTabJMAA2021} were derived from Navier-Stokes equations in  a three-dimensional thin domain, $\Omega^{\varepsilon}_t$, 
filled by a viscous fluid, that varies with time $t \in [0, T]$, given by
\begin{eqnarray}\Omega^{\varepsilon}_t&=&\left\{ (x_1^{\varepsilon},x_2^{\varepsilon},x_3^{\varepsilon})\in{R}^3:x_i(\xi_1,\xi_2,t)\leq x_i^{\varepsilon} \leq
x_i(\xi_1,\xi_2,t)+ h^\varepsilon(\xi_1,\xi_2,t)N_i(\xi_1,\xi_2,t),  \nonumber\right.\\
&&\left. (i=1,2,3), \ (\xi_1,\xi_2)\in D\subset \mathbb{R}^2
\right\} \label{eq-o-domain} \end{eqnarray} where
$\vec{X}_t(\xi_1,\xi_2)=\vec{X}(\xi_1,\xi_2,t)=(x_1(\xi_1,\xi_2,t),x_2(\xi_1,\xi_2,t),
x_3(\xi_1,\xi_2,t))$ is the lower
bound surface parametrization, $h^\varepsilon(\xi_1,\xi_2,t)$ is the gap between the two surfaces in
motion, and $\vec{N}(\xi_1,\xi_2,t)$ is the unit normal vector.

The lower bound surface is assumed to be regular
and the gap is assumed to be small with regard to the dimensions of the bound
surfaces. We take into account that the fluid film between the
surfaces is thin by introducing a small non-dimensional parameter 
$\varepsilon$, and setting that 
\begin{equation}
h^\varepsilon(\xi_1,\xi_2,t) = \varepsilon h(\xi_1,\xi_2,t), \quad
h(\xi_1,\xi_2,t) \ge h_0 > 0, \quad \forall \ (\xi_1,\xi_2)\in D\subset \mathbb{R}^2, \ \forall \ 
t\in [0,T].
\end{equation}

We introduce a reference domain 
\begin{equation}
\Omega=D \times [0,1] \label{eq-Omega}
\end{equation}
independent of
$\varepsilon$ and $t$, which is related to $\Omega^{\varepsilon}_t$ by the following change of variable: 
\begin{eqnarray}
    t^\varepsilon&=&t \label{eq-1-cv} \\
    x_i^\varepsilon &=& x_i(\xi_1,\xi_2,t)+\varepsilon \xi_3 h(\xi_1,\xi_2,t)N_i(\xi_1,\xi_2,t) \label{eq-2-cv}
\end{eqnarray}
where $(\xi_1,\xi_2)\in D$ and $\xi_3 \in[0,1]$.

We also define the basis $\left\{
\vec{a}_1,\vec{a}_2,\vec{a}_3\right\}$ 
\begin{eqnarray}
 \vec{a}_1(\xi_1,\xi_2,t)&=&\dfrac{\partial \vec{X}(\xi_1,\xi_2,t)}{\partial \xi_1} \label{base_a1} \\
 \vec{a}_2(\xi_1,\xi_2,t)&=&\dfrac{\partial \vec{X}(\xi_1,\xi_2,t)}{\partial \xi_2}\\
 \vec{a}_3(\xi_1,\xi_2,t)&=& \vec{N}(\xi_1,\xi_2,t) \label{base_a3}
\end{eqnarray}

The velocity, $\vec{u}^{\varepsilon}$, and  the external density of 
volume forces, $\vec{f}^\varepsilon$, can be written in the new basis \eqref{base_a1}-\eqref{base_a3} as follows, where we adopt the convention of summing over repeated indices from 1 to 3, except where otherwise indicated:
\begin{eqnarray}\vec{u}^{\varepsilon} &=& u_i^{\varepsilon}
\vec{e}_i = u_k(\varepsilon)\vec{a}_{k}, \quad 
u_i^{\varepsilon}= \left ( u_k(\varepsilon)\vec{a}_{k} \right ) \cdot \vec{e}_i = u_k(\varepsilon) a_{ki} \label{cambio_base_u}\\
\vec{f}^\varepsilon &=& f_i^{\varepsilon}\vec{e}_i = 
f_k(\varepsilon)\vec{a}_{k}, \quad f_i^{\varepsilon}= \left ( f_k(\varepsilon)\vec{a}_{k} \right ) \cdot \vec{e}_i = f_k(\varepsilon){a}_{ki}
\label{cambio_base_f}
\end{eqnarray}
where ${a}_{ki} = \vec{a}_{k} \cdot \vec{e}_i$.

Taking into account \eqref{eq-1-cv}-\eqref{cambio_base_f}, Navier-Stokes equations can be written in the reference domain $\Omega$ in the following way (in the next equation, repeated indices indicates summation from 1 to 3, except for index $n$, which take values from 1 to 2):
\begin{eqnarray}
&&\dfrac{ \partial u_k(\varepsilon)}{\partial t}  {a}_{ki} +
u_k(\varepsilon)\dfrac{
\partial {a}_{ki}}{\partial t} + \left({a}_{ki} \dfrac{
\partial u_k(\varepsilon)}{\partial \xi_n} + u_k(\varepsilon)\dfrac{ \partial {a}_{ki}}{\partial \xi_n} \right)\left[ -(\alpha_n \vec{a}_1 + \beta_n \vec{a}_2)\cdot\left(
\dfrac{\partial \vec{X}}{\partial t} + \varepsilon \xi_3 h
\dfrac{\partial \vec{a}_{3}}{\partial t} \right) \right]
\nonumber\\
&&\hspace*{+0.5cm}{}+\left({a}_{ki}\dfrac{
\partial u_k(\varepsilon)}{\partial \xi_3} + u_k(\varepsilon)\dfrac{
\partial {a}_{ki}}{\partial \xi_3} \right)\left(
-\dfrac{1}{\varepsilon h} \vec{a}_3 \cdot\dfrac{\partial
\vec{X}}{\partial t}- \dfrac{\xi_3}{ h} \dfrac{\partial h}{\partial
t}
\right)\nonumber\\
&&\hspace*{+0.5cm}{}+u_k(\varepsilon){a}_{kj} \left({a}_{qi} \dfrac{
\partial u_q(\varepsilon)}{\partial \xi_l} +u_q(\varepsilon)
\dfrac{
\partial {a}_{qi}}{\partial \xi_l} \right) \left(  \alpha_l
{a}_{1j} + \beta_l {a}_{2j}+ \gamma_l {a}_{3j}\right)\nonumber\\
&&\hspace*{+0.5cm}=-\dfrac{1}{\rho_0}\dfrac{
\partial p(\varepsilon)}{\partial \xi_l}\left(\alpha_l
{a}_{1i} + \beta_l {a}_{2i}+ \gamma_l {a}_{3i}\right) + \nu \left\{ \left[\dfrac{
\partial^2 (u_k(\varepsilon){a}_{ki})}{\partial \xi_l \partial \xi_m}
\left( \alpha_l {a}_{1j} + \beta_l {a}_{2j}+ \gamma_l {a}_{3j}
\right)\right.\right. \nonumber\\
&&\hspace*{+0.5cm}\left.\left.{}+ \dfrac{
\partial (u_k(\varepsilon){a}_{ki})}{\partial \xi_l}\dfrac{\partial}{\partial \xi_m}\left(
\alpha_l {a}_{1j} + \beta_l {a}_{2j}+ \gamma_l {a}_{3j} \right)
\right] \left( \alpha_m {a}_{1j} + \beta_m {a}_{2j}+ \gamma_m
{a}_{3j} \right) \right\} \nonumber\\
&&\hspace*{+0.5cm}{}+f_k(\varepsilon){a}_{ki}, \quad (i=1,2,3)\label{ec_ns_ij_alfa_beta}\\
&& \left({a}_{kj} \dfrac{
\partial u_k(\varepsilon)}{\partial \xi_l} + u_k(\varepsilon)\dfrac{
\partial {a}_{kj}}{\partial \xi_l}\right) \left(\alpha_l
{a}_{1j} + \beta_l {a}_{2j}+ \gamma_l {a}_{3j}\right) =0
\label{div_i_dr_alfa_beta}
\end{eqnarray}where $\alpha_l$, $\beta_l$ and $\gamma_l$ are defined in appendix \ref{ApendiceA} by expressions \eqref{alfaides}-\eqref{beta3n}.  We denote by $p$ the pressure, by $\rho_0$ the fluid density and by $\nu$ the kinematic viscosity. 

We begin assuming that $u_i(\varepsilon)$, $f_i(\varepsilon)$
($i=1,2,3$) and $p(\varepsilon)$ can be developed in powers of
$\varepsilon$, that is:
\begin{eqnarray}
&& u_i(\varepsilon) = u_i^0 + \varepsilon u_i^1 + \varepsilon^2
u_i^2 + \cdots \quad (i=1,2,3) \label{ansatz_1}\\
&&p(\varepsilon) =\varepsilon^{-2} p^{-2} + \varepsilon^{-1} p^{-1}
+p^0 + \varepsilon p^1 + \varepsilon^2 p^2 + \cdots \label{ansatz_2}\\
&& f_i(\varepsilon) = f_i^0 + \varepsilon f_i^1 + \varepsilon^2
f_i^2 + \cdots \quad (i=1,2,3) \label{ansatz_3}
\end{eqnarray}

As mentioned above, using asymptotic analysis we are able to derive two different models depending on the boundary conditions chosen.

In the  first place, if we assume that the fluid slips at the lower surface $(\xi_3=0)$, and at the upper surface $(\xi_3=1)$, but there is continuity in the normal direction, so the tangential velocities at the lower and upper surfaces are known, and the normal velocity of each of them must match the fluid velocity, we obtain
\begin{eqnarray}
&&\hspace*{-0.9cm} \dfrac{1}{\sqrt{A^0}} \textrm{div}\left(\dfrac{h^3 }{ \sqrt{A^0}} M \nabla
p^{-2} \right) =12\mu \dfrac{\partial h}{\partial t} +
12\mu\dfrac{h A^1}{A^0} \left(\dfrac{\partial \vec{X}}{\partial t}
\cdot
\vec{a}_3\right)\nonumber\\
&&\hspace*{-0.5cm}{} - 6\mu \nabla h\cdot (\vec{W}^{0}-\vec{V}^{0}) +
\dfrac{6\mu h}{\sqrt{A^0}} \textrm{div} (\sqrt{A^0}(\vec{W}^{0} + \vec{V}^{0})) \label{Reynolds_gen}
\end{eqnarray}
that can be considered a generalization of Reynolds equation. We denote by $V_1 \vec{a}_1+V_2 \vec{a}_2$
the tangential velocity
at the lower surface and, by $W_1\vec{a}_1+W_2\vec{a}_2$ the
tangential velocity at the upper surface, and we have
\begin{eqnarray}
\vec{V}(\varepsilon)&=&(V_1,V_2)=\vec{V}^0+O(\varepsilon)\\
\vec{W}(\varepsilon)&=&(W_1,W_2)=\vec{W}^0+O(\varepsilon)
\end{eqnarray}   

Coefficient $A^0$ and matrix $M$ are defined in appendix \ref{ApendiceA} (\eqref{A0}-\eqref{M}), and $\mu = \rho_0 \nu$ is the dynamic viscosity.

Once obtained $p^{-2}$ using \eqref{Reynolds_gen}, the following approximation of the three components of the velocity is yielded
\begin{eqnarray}
&&\hspace*{-0.2cm} u_i^0 =\dfrac{h^2 (\xi_3^2-\xi_3)}{2\mu} \sum_{k=1}^2J^{0,0}_{ik} \dfrac{
\partial p^{-2}}{\partial \xi_k}+\xi_3(W_i^{0}-
V_i^{0})+ V_i^{0}, \quad (i=1,2)\label{ui_0_lub}\\
&&\hspace*{-0.2cm}
u_3^0  =  \dfrac{\partial \vec{X}}{\partial t} \cdot \vec{a}_3
\label{u30_lub}
\end{eqnarray}
where $J^{0,0}_{ik}$ is given by \eqref{J}.

If instead of considering that the tangential and normal velocities are known on the upper and lower surfaces, we assume that the normal component of the traction on $\xi_3=0$ and on $\xi_3=1$ are known pressures (denoted by $\pi_0^\varepsilon$ and $\pi_1^\varepsilon$, respectively), and that the tangential component of the traction on these surfaces are friction forces depending on the value of the velocities on $\partial D$, then we obtain a shallow water model: 
\begin{eqnarray}
&& \hspace*{-0.5cm} u_i^0=W_i^{0}=V_i^{0} \quad (i=1,2) \label{u_i^0s}\\
&& \hspace*{-0.5cm} p^{-2}=p^{-1}=0\label{p-2p-1}\\
&& \hspace*{-0.5cm} p^0=\frac{2\mu}{h} \dfrac{\partial h}{\partial t} + \pi_0^0 \label{p0s-v2}
\\
&&\hspace*{-0.5cm}\dfrac{ \partial V_i^0}{\partial t}  + \sum_{l=1}^2 \left( V_l^0-C^0_l \right) \dfrac{
\partial V_i^0}{\partial \xi_l} +
\sum_{k=1}^2 
\left( R^0_{ik}+\sum_{l=1}^2
H^0_{ilk} V_l^0 \right) V_k^0=-\dfrac{1}{\rho_0} \sum_{l=1}^2 \dfrac{
\partial  \pi_0^0  }{\partial \xi_l}J^{0,0}_{il}  \nonumber \\ 
&&{} +\nu \left\{ \sum_{m=1}^2 \sum_{l=1}^2 
 \dfrac{
\partial^2 V_i^0 }{\partial \xi_m \partial \xi_l} J^{0,0}_{lm}
+ \sum_{k=1}^2 \sum_{l=1}^2  \dfrac{
\partial V_k^0 }{\partial \xi_l}( L^0_{kli} +\psi(h)^0_{ikl} )\right.\nonumber\\
&& \left. {}
+ \sum_{k=1}^2 V_k^0 ({P}_{ik}^0+\chi(h)^0_{ik})  + {\kappa}(h)^0_i  \right\}+ {F}^0_i(h)-
 Q^0_{i3} \left ( \frac{\partial \vec{X}}{\partial t}
\cdot \vec{a}_3\right )  \quad (i=1,2)\label{ec_Vi0-v2} \end{eqnarray}
where coefficients $C^0_l$, $R^0_{ik}$, $H^0_{ilk}$, $J^{0,0}_{lm}$, $L^0_{kli}$, $\psi(h)^0_{ikl}$, $P^0_{ik}$, $\chi(h)^0_{ik}$, ${\kappa}(h)^0_i$, $F^0_i(h)$ and $Q^0_{i3}$ are defined in appendix \ref{ApendiceA} (in \eqref{C}, \eqref{R}, \eqref{H}, \eqref{J}, \eqref{L}, \eqref{psi}, \eqref{P}, \eqref{chiik}, \eqref{kappa}, \eqref{F} and \eqref{Q3} respectively),  and $\pi_0^0$ is the term of order zero on $\varepsilon$ of $\pi_0^\varepsilon$, that is,  $\pi_0^\varepsilon = \pi_0^0 + O(\varepsilon)$.

\section{New hypothesis about the dependence of the solution on $\xi_3$} \label{NewAssumptions}

If we carefully observe the steps of the proofs in the previous work \cite{RodTabJMAA2021}, we can see that $p^k\ (k=-2, -1, 0, 1)$, and $u_i^k\ (i=1,2,3; k=0,1)$ are polynomials in $\xi_3$ of at most degree three. Because of this, we are going to assume that, for $\varepsilon$ small enough, the following equalities are true:

\begin{eqnarray}
u_i(t,\xi_1,\xi_2,\xi_3,\varepsilon)&=&\sum_{n=0}^3 \xi_3^n \bar{u}_i^n(t,\xi_1,\xi_2,\varepsilon), \quad (i=1,2,3) \label{ui_pol_xi3}\\
p(t,\xi_1,\xi_2,\xi_3,\varepsilon)&=&\sum_{n=0}^3 \xi_3^n \bar{p}^n(t,\xi_1,\xi_2,\varepsilon),\label{p_pol_xi3} \\
{f_i}(t,\xi_1,\xi_2,\xi_3,\varepsilon)&=&\sum_{n=0}^{\infty} \xi_3^n \bar{f}_i^n(t,\xi_1,\xi_2,\varepsilon).\label{f_pol_xi3} 
\end{eqnarray}

We want to point out that the previous hypothesis is equivalent to neglecting in \eqref{ansatz_1}-\eqref{ansatz_2} the terms in $O(\varepsilon^2)$ when $\varepsilon$ is small.

Using expressions \eqref{ui_pol_xi3}-\eqref{f_pol_xi3} and \eqref{alfaides}-\eqref{gammades} we can rewrite equations \eqref{ec_ns_ij_alfa_beta}-\eqref{div_i_dr_alfa_beta} as follows (repeated indices $k$ and $j$ indicate summation from 1 to 3, while repeated indices $l$ and $m$ indicate summation from 1 to 2):
\begin{eqnarray}
&&\hspace*{-0.5cm} \sum_{n=0}^3 \xi_3^n \dfrac{ \partial \bar{u}_k^n}{\partial t} \vec{a}_{k} +
\sum_{n=0}^3 \xi_3^n \bar{u}_k^n \dfrac{
\partial \vec{a}_{k}}{\partial t} \nonumber\\
&&{}- \left(\vec{a}_{k}\sum_{n=0}^3 \xi_3^n  \dfrac{
\partial \bar{u}_k^n}{\partial \xi_l} + \sum_{n=0}^3 \xi_3^n \bar{u}_k^n\dfrac{ \partial \vec{a}_{k}}{\partial \xi_l} \right)\left[ \sum_{r=0}^{\infty} (\varepsilon \xi_3 h)^r \left(\alpha_l^r \vec{a}_1 +  \beta_l^r \vec{a}_2\right) \cdot\left(
\dfrac{\partial \vec{X}}{\partial t} + \varepsilon \xi_3 h
\dfrac{\partial \vec{a}_{3}}{\partial t} \right) \right]
\nonumber\\
&&{}-\vec{a}_{k} \sum_{n=1}^3 n \xi_3^{n-1} \bar{u}_k^n\left[\dfrac{1}{\varepsilon h} \left(\vec{a}_3 \cdot\dfrac{\partial
\vec{X}}{\partial t} \right)+ \dfrac{\xi_3}{ h} \dfrac{\partial h}{\partial
t}\right.\nonumber\\
&&\left.{}+ \sum_{r=0}^{\infty} \varepsilon^{r} \xi_3^{r+1} h^{r-1} \left( \alpha_3^{r} \vec{a}_1 + \beta_3^{r} \vec{a}_2 \right) \cdot\left(
\dfrac{\partial \vec{X}}{\partial t} + \varepsilon \xi_3 h
\dfrac{\partial \vec{a}_{3}}{\partial t} \right)
\right]\nonumber\\
&&{}+\sum_{n=0}^3 \xi_3^n \bar{u}_k^n  \left(\vec{a}_{q}\sum_{d=0}^3 \xi_3^{d}  \dfrac{
\partial \bar{u}_q^{d}}{\partial \xi_l} + \sum_{d=0}^3 \xi_3^{d} \bar{u}_q^{d}\dfrac{ \partial \vec{a}_{q}}{\partial \xi_l} \right)   \sum_{r=0}^{\infty} (\varepsilon \xi_3 h)^{r} \left(\alpha_l^{r}
(\vec{a}_k \cdot \vec{a}_{1}) + \beta_l^{r} (\vec{a}_k \cdot \vec{a}_{2})\right))\nonumber\\
&&{}+\sum_{n=0}^3 \xi_3^n \bar{u}_k^n \left(\vec{a}_{q} \sum_{d=1}^3 d \xi_3^{d-1} \bar{u}_q^{d} \right) \left(  \sum_{r=0}^{\infty} \varepsilon^n \xi_3^{r+1} h^{r-1} \left( \alpha_3^{r}
(\vec{a}_k \cdot \vec{a}_{1}) + \beta_3^{r} (\vec{a}_k \cdot \vec{a}_{2})\right)\right.\nonumber\\
&&\left.{}+ \dfrac{1}{\varepsilon h} (\vec{a}_k \cdot \vec{a}_{3})\right)=-\dfrac{1}{\rho_0} \sum_{n=0}^3 \xi_3^n \dfrac{
\partial  \bar{p}^n}{\partial \xi_l} \sum_{r=0}^{\infty} (\varepsilon \xi_3 h)^{r} \left(\alpha_l^{r} \vec{a}_1 +  \beta_l^{r} \vec{a}_2\right) \nonumber\\
&&{} -\dfrac{1}{\rho_0} \sum_{n=1}^3 n \xi_3^{n-1} \bar{p}^n\left( \sum_{r=0}^{\infty} \varepsilon^{r} \xi_3^{{r}+1} h^{r-1} \left( \alpha_3^{r} \vec{a}_1 + \beta_3^{r} \vec{a}_2 \right) + \dfrac{1}{\varepsilon h} \vec{a}_{3}\right)\nonumber\\
&&{} + \nu \left[ \displaystyle\sum_{n=0}^3 \xi_3^n\dfrac{
\partial^2 (\vec{a}_{k}  \bar{u}_k^n)}{\partial \xi_l \partial \xi_m}
\sum_{r=0}^{\infty} (\varepsilon \xi_3 h)^{r} \left( \alpha_l^{r} {a}_{1j} + \beta_l^{r} {a}_{2j}
\right)\right. \nonumber\\
&&\left.{}+ 2\sum_{n=1}^3 n\xi_3^{n-1}\dfrac{
\partial (\vec{a}_{k}  \bar{u}_k^n)}{ \partial \xi_m}
\left(\sum_{r=0}^{\infty} \varepsilon^{r}
\xi_3^{r+1} h^{r-1} \left( \alpha_3^{r} {a}_{1j} + \beta_3^{r}
{a}_{2j}\right)+ \dfrac{1}{\varepsilon h}
{a}_{3j}
\right)\right. \nonumber\\
&&\left.{}+ \displaystyle\sum_{n=0}^3 \xi_3^n\dfrac{
\partial (\vec{a}_{k}  \bar{u}_k^n)}{\partial \xi_l}\dfrac{\partial}{\partial \xi_m} \left(\sum_{r=0}^{\infty} (\varepsilon \xi_3 h)^{r} \left( \alpha_l^{r} {a}_{1j} + \beta_l^{r} {a}_{2j}
\right) \right)\right. \nonumber\\
&&\left.{}+  \vec{a}_{k}
\displaystyle\sum_{n=1}^3 n \xi_3^{n-1}
\bar{u}_k^n\dfrac{\partial}{\partial \xi_m}\left(\sum_{r=0}^{\infty} \varepsilon^{r}
\xi_3^{r+1} h^{r-1} \left( \alpha_3^{r} {a}_{1j} + \beta_3^{r}
{a}_{2j}\right)+ \dfrac{1}{\varepsilon h}
{a}_{3j}  \right)
\right]\nonumber\\
&& \cdot \sum_{s=0}^{\infty} (\varepsilon \xi_3 h)^{s} \left( \alpha_m^{s} {a}_{1j} + \beta_m^{s} {a}_{2j}
\right) \nonumber\\
&&{}+ \nu \left[\vec{a}_{k} \sum_{n=2}^3 n(n-1) \xi_3^{n-2}
\bar{u}_k^n \left( \sum_{r=0}^{\infty} \varepsilon^{r}
\xi_3^{r+1} h^{r-1} \left( \alpha_3^{r} {a}_{1j} + \beta_3^{r}
{a}_{2j}\right)+ \dfrac{1}{\varepsilon h}
{a}_{3j}
\right)\right. \nonumber\\
&&\left.{}+ \displaystyle\sum_{n=0}^3 \xi_3^n \dfrac{
\partial (\vec{a}_{k} \bar{u}_k^n)}{\partial \xi_l} 
\sum_{r=1}^{\infty} r \varepsilon^r \xi_3^{r-1} h^{r} \left( \alpha_l^{r} {a}_{1j} + \beta_l^{r} {a}_{2j}
\right) \right. \nonumber\\
&&\left.{}+\vec{a}_{k} \displaystyle\sum_{n=1}^3 n \xi_3^{n-1}
\bar{u}_k^n \left( \sum_{r=0}^{\infty} (r+1)\varepsilon^{r} \xi_3^{r}
h^{r-1} \left( \alpha_3^{r} {a}_{1j} + \beta_3^{r} {a}_{2j}
\right) \right)\right]\nonumber\\
&& \cdot\left( \sum_{s=0}^{\infty} \varepsilon^{s}
\xi_3^{s+1} h^{s-1} \left( \alpha_3^{s} {a}_{1j} + \beta_3^{s}
{a}_{2j}\right) + \dfrac{1}{\varepsilon h}
{a}_{3j} \right)+
\sum_{n=0}^{\infty} \xi_3^n \bar{f}_k^n \vec{a}_{k}\label{ec_ns_ij_alfa_beta_pol_xi3}
\\
&&\hspace*{-0.5cm}\sum_{n=0}^3 \xi_3^n \dfrac{
\partial \bar{u}_k^n}{\partial \xi_l} \sum_{r=0}^{\infty} (\varepsilon \xi_3 h)^{r} \left( \alpha_l^{r}
(\vec{a}_k \cdot \vec{a}_{1}) + \beta_l^{r} (\vec{a}_k \cdot \vec{a}_{2})\right) \nonumber\\
&&{} +\sum_{n=0}^3 \xi_3^n \bar{u}_k^n\dfrac{
\partial {a}_{ki}}{\partial \xi_l} \sum_{r=0}^{\infty} (\varepsilon \xi_3 h)^{r} \left( \alpha_l^{r} {a}_{1i} + \beta_l^{r} {a}_{2i}
\right) \nonumber\\
&&{}+ \sum_{n=1}^3 n \xi_3^{n-1 } \bar{u}_k^n \left(  \sum_{r=0}^{\infty} \varepsilon^{r}
\xi_3^{r+1} h^{r-1} \left( \alpha_3^{r}
(\vec{a}_k \cdot \vec{a}_{1}) + \beta_3^{r} (\vec{a}_k \cdot \vec{a}_{2})\right) + \dfrac{1}{\varepsilon h} (\vec{a}_k \cdot \vec{a}_{3})\right) =0
\label{div_i_dr_alfa_beta_pol_xi3}
\end{eqnarray}
and identify the terms multiplied by the different powers of $\xi_3$.

Let us begin identifying the terms multiplied by $\xi_3^0$ in equation \eqref{div_i_dr_alfa_beta_pol_xi3}:
\begin{eqnarray}
\bar{u}_3^1&=&-\varepsilon h\left[\dfrac{
\partial \bar{u}_1^0}{\partial \xi_1} + \dfrac{
\partial \bar{u}_2^0}{\partial \xi_2} + \dfrac{1}{2A^0} \left(\dfrac{
\partial A^0}{\partial \xi_1} \bar{u}_1^0 +\dfrac{
\partial A^0}{\partial \xi_2} \bar{u}_2^0  \right) {+\bar{u}_3^0 \dfrac{A^1}{A^0}} \right]
\label{div_xi3_0_res}\\
&=&-\varepsilon h\left[ \textrm{div}(\vec{u}^0)+ \dfrac{1}{2A^0} \nabla A^0 \cdot \vec{u}^0 {+\bar{u}_3^0  \dfrac{A^1}{A^0}}\right] =-\dfrac{\varepsilon h}{\sqrt{A^0}} \textrm{div}(\sqrt{A^0}\vec{u}^0)  {-\varepsilon h \bar{u}_3^0  \dfrac{A^1}{A^0}}\label{u_3^1}
\end{eqnarray}where
$\vec{u}^0=(\bar{u}_1^0,\bar{u}_2^0)$, $A^0$ and $A^1$ are given by \eqref{A0} and \eqref{A^1}.

Now, we identify the terms multiplied by $\xi_3^0$ in \eqref{ec_ns_ij_alfa_beta_pol_xi3}.  In equations \eqref{ec_ns_xi3_0}-\eqref{ec_ns_xi3_0_a3}, below,  repeated indices $k$ and $j$ indicate again summation from 1 to 3, while repeated indices $l$ and $m$ indicate summation from 1 to 2. 
\begin{eqnarray}
&&\hspace*{-0.5cm} \dfrac{ \partial \bar{u}_k^0}{\partial t} \vec{a}_{k} +
 \bar{u}_k^0 \dfrac{
\partial \vec{a}_{k}}{\partial t} - \left(\vec{a}_{k}  \dfrac{
\partial \bar{u}_k^0}{\partial \xi_l} +  \bar{u}_k^0\dfrac{ \partial \vec{a}_{k}}{\partial \xi_l} \right) C^0_l-\vec{a}_{k}   \bar{u}_k^1\dfrac{1}{\varepsilon h} \left(\vec{a}_3 \cdot\dfrac{\partial
\vec{X}}{\partial t} \right)
\nonumber\\
&&{}
+ \bar{u}_k^0 \left(\vec{a}_{k} \dfrac{
\partial \bar{u}_k^0}{\partial \xi_l} +  \bar{u}_k^0\dfrac{ \partial \vec{a}_{k}}{\partial \xi_l} \right)     \left(\alpha_l^0
(\vec{a}_k \cdot \vec{a}_{1}) + \beta_l^0 (\vec{a}_k \cdot \vec{a}_{2})\right))+ \bar{u}_k^0 \left(\vec{a}_{k}  \bar{u}_k^1\right) \left(  \dfrac{1}{\varepsilon h} (\vec{a}_k \cdot \vec{a}_{3})\right)\nonumber\\
&&{}=-\dfrac{1}{\rho_0}  \dfrac{
\partial  \bar{p}^0}{\partial \xi_l}   \left(\alpha_l^0\vec{a}_1 +  \beta_l^0 \vec{a}_2\right) -\dfrac{1}{\rho_0} \dfrac{1}{\varepsilon h} \bar{p}^1 \vec{a}_{3} + \nu \left[\dfrac{
\partial^2 (\vec{a}_{k}  \bar{u}_k^0)}{\partial \xi_l \partial \xi_m}
 \left( \alpha_l^0 {a}_{1j} + \beta_l^0 {a}_{2j}
\right)+\dfrac {2}{\varepsilon h} \dfrac{
\partial (\vec{a}_{k}  \bar{u}_k^1)}{ \partial \xi_m}
{a}_{3j}\right. \nonumber\\
&&\left.{}
+ \dfrac{
\partial (\vec{a}_{k}  \bar{u}_k^0)}{\partial \xi_l}\dfrac{\partial}{\partial \xi_m} \left( \alpha_l^0 {a}_{1j} + \beta_l^0 {a}_{2j}
\right)+  \vec{a}_{k}
\bar{u}_k^1\dfrac{\partial}{\partial \xi_m}\left( \dfrac{1}{\varepsilon h}
{a}_{3j}  \right)
\right]  \left( \alpha_m^0 {a}_{1j} + \beta_m^0 {a}_{2j}
\right) \nonumber\\
&&{}+  \dfrac{\nu}{\varepsilon h}
{a}_{3j} \left[  \dfrac{2}{\varepsilon h}\vec{a}_{k}
\bar{u}_k^2
{a}_{3j}
+  \dfrac{
\partial (\vec{a}_{k} \bar{u}_k^0)}{\partial \xi_l}
 \varepsilon  h \left( \alpha_l^1 {a}_{1j} + \beta_l^1 {a}_{2j}
\right)+\vec{a}_{k}
\bar{u}_k^1 \left(
h^{-1} \left( \alpha_3^0 {a}_{1j} + \beta_3^0 {a}_{2j}
\right)\right)\right]+ {\bar{f}_k^0}\vec{a}_{k}\nonumber\\\label{ec_ns_xi3_0}
\end{eqnarray}
where $C^0_l$ is given by \eqref{C}.

We multiply \eqref{ec_ns_xi3_0} by $\vec{a}_i$ ($i=1,2,3$), and we yield:
\begin{eqnarray}
&&\hspace*{-0.5cm} \dfrac{ \partial \bar{u}_1^0}{\partial t} (\vec{a}_1 \cdot \vec{a}_i)+ \dfrac{ \partial \bar{u}_2^0}{\partial t}  (\vec{a}_2 \cdot \vec{a}_i) +
 \bar{u}_1^0 \left(\vec{a}_i \cdot \dfrac{
\partial \vec{a}_{1}}{\partial t}\right) +
 \bar{u}_2^0 \left(\vec{a}_i \cdot \dfrac{
\partial \vec{a}_{2}}{\partial t}\right)+
 \bar{u}_3^0 \left(\vec{a}_i \cdot \dfrac{
\partial \vec{a}_{3}}{\partial t}\right) \nonumber\\
&&{}- \left( \dfrac{ \partial \bar{u}_1^0}{\partial \xi_l}  (\vec{a}_1 \cdot \vec{a}_i) + \dfrac{ \partial \bar{u}_2^0}{\partial \xi_l}  (\vec{a}_2 \cdot \vec{a}_i) + \sum_{k=1}^3
 \bar{u}_k^0 \left(\vec{a}_i \cdot \dfrac{
\partial \vec{a}_{k}}{\partial \xi_l}\right)\right)
C^0_l
\nonumber\\
&&{}+ \dfrac{1}{\varepsilon h} \left[\bar{u}_3^0 -   \left(\vec{a}_3 \cdot\dfrac{\partial
\vec{X}}{\partial t} \right)\right]\left(  \bar{u}_1^1 (\vec{a}_1 \cdot \vec{a}_i) +  \bar{u}_2^1 (\vec{a}_2 \cdot \vec{a}_i) \right) \nonumber\\
&&{} + \bar{u}_1^0 \left( \dfrac{
\partial \bar{u}_1^0}{\partial \xi_1}(\vec{a}_1 \cdot \vec{a}_i)  +  \dfrac{
\partial \bar{u}_2^0}{\partial \xi_1} (\vec{a}_1 \cdot \vec{a}_i) + \sum_{k=1}^3
 \bar{u}_k^0 \left(\vec{a}_i \cdot \dfrac{
\partial \vec{a}_{k}}{\partial \xi_1}\right) \right)  \nonumber\\
&&{}
+ \bar{u}_2^0 \left(\dfrac{
\partial \bar{u}_1^0}{\partial \xi_2} (\vec{a}_1 \cdot \vec{a}_i)  + \dfrac{
\partial \bar{u}_2^0}{\partial \xi_2} (\vec{a}_1 \cdot \vec{a}_i)  +\sum_{k=1}^3
 \bar{u}_k^0 \left(\vec{a}_i \cdot \dfrac{
\partial \vec{a}_{k}}{\partial \xi_2}\right)  \right)   \nonumber\\
&&{}
 =-\dfrac{1}{\rho_0}  \dfrac{
\partial  \bar{p}^0}{\partial \xi_i}  + \nu \vec{a}_i\cdot   \left[\dfrac{
\partial^2 (\vec{a}_{k}  \bar{u}_k^0)}{\partial \xi_l \partial \xi_m}
 \left( \alpha_l^0 {a}_{1j} + \beta_l^0 {a}_{2j}
\right)
\right. \nonumber\\
&&\left.{}+ \dfrac{
\partial (\vec{a}_{k}  \bar{u}_k^0)}{\partial \xi_l}\dfrac{\partial}{\partial \xi_m} \left( \alpha_l^0 {a}_{1j} + \beta_l^0 {a}_{2j}
\right)+  \vec{a}_{k}
\bar{u}_k^1\dfrac{\partial}{\partial \xi_m}\left( \dfrac{1}{\varepsilon h}
{a}_{3j}  \right)
\right]  \left( \alpha_m^0 {a}_{1j} + \beta_m^0 {a}_{2j}
\right)  \nonumber\\
&&{}+  \dfrac{2\nu}{\varepsilon^2 h^2}
 \left(  \bar{u}_1^2(\vec{a}_1 \cdot \vec{a}_i)  + \bar{u}_2^2 (\vec{a}_1 \cdot \vec{a}_i) 
\right)+ (\vec{a}_1 \cdot \vec{a}_i) \bar{f}_1^0+ (\vec{a}_1 \cdot \vec{a}_i) {\bar{f}_2^0} \quad (i=1,2)\label{ec_ns_xi3_0_ai}
\\
&&\hspace*{-0.5cm} \dfrac{ \partial \bar{u}_3^0}{\partial t}  +
 \bar{u}_1^0 \left(\vec{a}_3 \cdot \dfrac{
\partial \vec{a}_{1}}{\partial t}\right) +
 \bar{u}_2^0 \left(\vec{a}_3 \cdot \dfrac{
\partial \vec{a}_{2}}{\partial t}\right)
\nonumber\\
&&{} -  \left(  \dfrac{
\partial \bar{u}_3^0}{\partial \xi_l} +
 \bar{u}_1^0 \left(\vec{a}_3 \cdot \dfrac{
\partial \vec{a}_{1}}{\partial \xi_l}\right) +
 \bar{u}_2^0 \left(\vec{a}_3 \cdot \dfrac{
\partial \vec{a}_{2}}{\partial \xi_l} \right)  \right) C^0_l
\nonumber\\
&&{}+   \dfrac{1}{\varepsilon h} \bar{u}_3^1 \left[\bar{u}_3^0 -\left(\vec{a}_3 \cdot\dfrac{\partial
\vec{X}}{\partial t} \right)\right] + \bar{u}_1^0 \left[\dfrac{
\partial \bar{u}_3^0}{\partial \xi_1} +  \bar{u}_1^0 \left(\vec{a}_3 \cdot \dfrac{ \partial \vec{a}_{1}}{\partial \xi_1}\right)+  \bar{u}_2^0 \left(\vec{a}_3 \cdot \dfrac{ \partial \vec{a}_{2}}{\partial \xi_1}\right) \right]
 \nonumber\\
&&{}+ \bar{u}_2^0 \left[\dfrac{
\partial \bar{u}_3^0}{\partial \xi_2} +  \bar{u}_1^0 \left(\vec{a}_3 \cdot \dfrac{ \partial \vec{a}_{1}}{\partial \xi_2}\right)+  \bar{u}_2^0 \left(\vec{a}_3 \cdot \dfrac{ \partial \vec{a}_{2}}{\partial \xi_2}\right) \right]
  = -\dfrac{1}{\rho_0} \dfrac{1}{\varepsilon h} \bar{p}^1 \nonumber\\
&&{} + \nu \vec{a}_3\cdot  \left[\dfrac{
\partial^2 (\vec{a}_{k}  \bar{u}_k^0)}{\partial \xi_l \partial \xi_m}
 \left( \alpha_l^0 {a}_{1j} + \beta_l^0 {a}_{2j}
\right)
+ \dfrac{
\partial (\vec{a}_{k}  \bar{u}_k^0)}{\partial \xi_l}\dfrac{\partial}{\partial \xi_m} \left( \alpha_l^0 {a}_{1j} + \beta_l^0 {a}_{2j}
\right)\right. \nonumber\\
&&\left.{}+  \vec{a}_{k}
\bar{u}_k^1\dfrac{\partial}{\partial \xi_m}\left( \dfrac{1}{\varepsilon h}
{a}_{3j}  \right)
\right]  \left( \alpha_m^0 {a}_{1j} + \beta_m^0 {a}_{2j}
\right)+  \dfrac{2\nu}{\varepsilon^2 h^2}
\bar{u}_3^2
+ {\bar{f}_3^0}\label{ec_ns_xi3_0_a3}
\end{eqnarray}

Next, we multiply equation \eqref{ec_ns_xi3_0_ai} ($i=1$) by $\alpha^0_k$ and we add equation \eqref{ec_ns_xi3_0_ai} ($i=2$) multiplied by $\beta^0_k$ for $k=1,2$ to get these two equations:
\begin{eqnarray}
&&\hspace*{-0.5cm} \dfrac{ \partial \bar{u}_i^0}{\partial t}+\sum_{l=1}^2  \dfrac{ \partial \bar{u}_i^0}{\partial \xi_l} (\bar{u}_l^0 -C^0_l)+ \sum_{k=1}^3
 \bar{u}_k^0  \left(  Q^0_{ik} + \sum_{l=1}^2 \bar{u}_l^0 H^0_{ilk} \right)
+ \dfrac{1}{\varepsilon h}  \bar{u}_i^1  \left[\bar{u}_3^0 - \left(\vec{a}_3 \cdot\dfrac{\partial
\vec{X}}{\partial t} \right)\right] \nonumber\\
&&{}
 =-\dfrac{1}{\rho_0}  \left(\alpha^0_i\dfrac{
\partial  \bar{p}^0}{\partial \xi_1} +\beta^0_i \dfrac{
\partial  \bar{p}^0}{\partial \xi_2}\right)  + \nu \left[ \sum_{m=1}^2 \sum_{l=1}^2 \dfrac{
\partial^2 \bar{u}_i^0}{\partial \xi_l \partial \xi_m} J^{0,0}_{lm} \right. \nonumber\\
&&\left.{}+ \sum_{k=1}^3 \sum_{l=1}^2 \dfrac{
\partial \bar{u}_k^0}{ \partial \xi_l} L^0_{kli}
+  \sum_{k=1}^3  \bar{u}^0_k {S}^0_{ik}+  \dfrac{1}{\varepsilon  h} \dfrac{A^1}{A^0}
\bar{u}_i^1  +  \dfrac{2}{\varepsilon^2 h^2}
 \bar{u}_i^2
\right]   + {\bar{f}_i^0} \quad (i=1,2)\label{ec_ui0_coefs}\end{eqnarray}
where coefficients $H^0_{ilk}$, $J^{0,0}_{lm}$, $L^0_{kli}$, $Q^0_{ik}$ and $S^0_{ik}$ are given by \eqref{H}, \eqref{J}, \eqref{L},  \eqref{Q} and \eqref{S} respectively.

Equation \eqref{ec_ns_xi3_0_a3} could be written in a more compact form using the coefficients $L^0_{kl3}$, $Q^0_{3k}$ and $S^0_{3k}$  defined in \eqref{L3}, \eqref{Q3} and \eqref{S3}:
\begin{eqnarray}
&&\hspace*{-0.5cm} \dfrac{ \partial \bar{u}_3^0}{\partial t} +\sum_{l=1}^2  \dfrac{ \partial \bar{u}_3^0}{\partial \xi_l} (\bar{u}_l^0 -C^0_l)
 + \sum_{k=1}^2
 \bar{u}_k^0\left[ Q^0_{3k}+ \sum_{l=1}^2 \bar{u}_l^0 \left(\vec{a}_3 \cdot \dfrac{ \partial \vec{a}_{l}}{\partial \xi_k}\right) \right]
\nonumber\\
&&{}+ \dfrac{1}{\varepsilon h}  \bar{u}_3^1  \left[\bar{u}_3^0 - \left(\vec{a}_3 \cdot\dfrac{\partial
\vec{X}}{\partial t} \right)\right] = -\dfrac{1}{\rho_0} \dfrac{1}{\varepsilon h} \bar{p}^1  + \nu \left[ \sum_{m=1}^2 \sum_{l=1}^2\dfrac{
\partial^2   \bar{u}_3^0}{\partial \xi_l \partial \xi_m}
J^{0,0}_{lm} \right. \nonumber\\
&&\left.{} +\sum_{k=1}^3 \sum_{l=1}^2 \dfrac{
\partial  \bar{u}_k^0}{\partial \xi_l } L^0_{kl3} + \sum_{k=1}^3  \bar{u}_k^0 {{S}^0_{3k}}+ \dfrac{1}{\varepsilon h} \dfrac{A^1}{A^0}
\bar{u}_3^1 +  \dfrac{2}{\varepsilon^2 h^2}
\bar{u}_3^2
\right]
+ {\bar{f}_3^0}\label{ec_ns_u30_Res}
\end{eqnarray}   

Subsequently, we identify the terms multiplied by $\xi_3$ in \eqref{ec_ns_ij_alfa_beta_pol_xi3}-\eqref{div_i_dr_alfa_beta_pol_xi3} and, following the steps \eqref{ec_ns_xi3_0_ai}-\eqref{ec_ns_u30_Res} we yield:
\begin{eqnarray}
&&\hspace*{-0.5cm}\bar{u}_3^2=- \dfrac{\varepsilon h}{2\sqrt{A^0}} \textrm{div}(\sqrt{A^0}\vec{u}^1) + \dfrac{\varepsilon }{2} \nabla h \cdot \vec{u}^1  -\varepsilon h \bar{u}_3^1 \dfrac{A^1}{A^0} -\dfrac{\varepsilon^2 h^2}{2}  \sum_{l=1}^2  \sum_{k=1}^3 \left( \dfrac{
\partial \bar{u}_k^0}{\partial \xi_l}  B^1_{lk}+   \bar{u}_k^0  H^1_{llk}\right)
\label{u32}\\
&&\hspace*{-0.5cm}  \dfrac{ \partial \bar{u}_i^1}{\partial t}+ \sum_{l=1}^2  \dfrac{ \partial \bar{u}_i^1}{\partial \xi_l} \left(  \bar{u}_l^0    -  C_l^0 \right)   + \sum_{k=1}^3 \bar{u}_k^1\left[ Q^0_{ik}+\sum_{l=1}^2  \bar{u}_l^0\left(H^0_{ilk}  -\dfrac{\delta_{ki}}{ h}   \dfrac{\partial h}{\partial \xi_l} \right) \right]\nonumber\\
&&{}+\sum_{l=1}^2 \bar{u}_l^1 \left( \dfrac{
\partial \bar{u}_i^0}{\partial \xi_l}  + \sum_{k=1}^3 \bar{u}_k^0 H^0_{ilk}
 \right) -\dfrac{1}{ h}  \bar{u}_i^1 \left( \dfrac{\partial
h}{\partial t}+C_3^0\right)\nonumber\\
&&{}+ \varepsilon  h \sum_{l=1}^2\left(  \dfrac{ \partial \bar{u}_i^0}{\partial
\xi_l}    +\sum_{k=1}^3
\bar{u}_k^0 H^0_{ilk}\right)\left[  \sum_{m=1}^2 B^1_{lm} \bar{u}_m^0 - C^{1,0}_l\right]
\nonumber\\
&&{}+\dfrac{2}{\varepsilon h}\bar{u}_i^2 \left[\bar{u}_3^0 -\left(\vec{a}_3 \cdot\dfrac{\partial
\vec{X}}{\partial t} \right)   \right]  +
\dfrac{1}{\varepsilon h}  \bar{u}_i^1 \bar{u}_3^1\nonumber\\
&&{}=-\dfrac{1}{\rho_0} \sum_{l=1}^2  \dfrac{
\partial  \bar{p}^1}{\partial \xi_l} J^{0,0}_{il}   - \dfrac{\varepsilon  h}{\rho_0} \sum_{l=1}^2  \dfrac{
\partial  \bar{p}^0}{\partial \xi_l} J^{0,1}_{il} -\dfrac{1}{\rho_0 h}\bar{p}^1 J^{0,0}_{i3}\nonumber\\
&&{} + \nu \left\{ \sum_{m=1}^2\sum_{l=1}^2 \dfrac{
\partial^2 \bar{u}_i^1}{\partial \xi_l \partial \xi_m}
J^{0,0}_{lm} + \sum_{l=1}^2  \sum_{k=1}^3  \dfrac{
\partial  \bar{u}_k^1}{\partial \xi_l } L^{1,0}_{kli}(h)+\sum_{k=1}^3   \bar{u}_k^1{S}^{1,0}_{ik} (h) \right.\nonumber\\
&&\left.{}+  \varepsilon \left[ 2 h \sum_{m=1}^2 \sum_{l=1}^2 \dfrac{
\partial^2  \bar{u}_i^0}{\partial \xi_l \partial \xi_m} J^{1,0}_{lm} + 
\sum_{l=1}^2 \sum_{k=1}^3  \dfrac{
\partial \bar{u}_k^0}{\partial \xi_l}L^{0,1}_{kli}(h)  + \sum_{k=1}^3  \bar{u}_k^0 S_{ik}^{0,1}(h)\right]\right.\nonumber\\
&&\left.{}+ \dfrac{2}{\varepsilon h}  \bar{u}_i^2 \dfrac{A^1}{A^0} + \dfrac{6}{\varepsilon^2 h^2}
\bar{u}_i^3 \right\}+ \bar{f}_i^1\quad (i=1,2)\label{ec_u1i}\\
&&\hspace*{-0.5cm}  \dfrac{ \partial \bar{u}_3^1}{\partial t}  + \sum_{l=1}^2   \dfrac{
\partial \bar{u}_3^1}{\partial \xi_l}\left(\bar{u}_l^0 -C_l^0\right) +
\sum_{k=1}^2 \bar{u}_k^1 \left[Q^0_{3k} + \sum_{l=1}^2  \left(\vec{a}_3 \cdot\dfrac{
\partial \vec{a}_{k}}{\partial \xi_l} \right)  \bar{u}_l^0 \right]\nonumber\\
&&{}+\sum_{l=1}^2 \bar{u}_l^1 \left(  \dfrac{
\partial \bar{u}_3^0}{\partial \xi_l} + \sum_{k=1}^2 \bar{u}_k^0 \left(\vec{a}_3 \cdot \dfrac{ \partial \vec{a}_{k}}{\partial \xi_l}\right)
\right) \nonumber\\
&&{}+ \varepsilon  h  \sum_{l=1}^2 \left[ \dfrac{
\partial \bar{u}_3^0}{\partial \xi_l} +
\sum_{k=1}^2 \bar{u}_k^0 \left(\vec{a}_3 \cdot\dfrac{
\partial \vec{a}_{k}}{\partial \xi_l} \right)  \right] \left(\sum_{m=1}^2 B^1_{lm} \bar{u}_m^0- C^{1,0}_l \right)\nonumber\\
&&{}+ \dfrac{2}{\varepsilon h}    \bar{u}_3^2\left[ \bar{u}_3^0 -\left(\vec{a}_3
\cdot\dfrac{\partial \vec{X}}{\partial t} \right) \right]+  \dfrac{1}{\varepsilon h}
\bar{u}_3^1 \left[\bar{u}_3^1-\varepsilon\left(\dfrac{\partial h}{\partial t}+ C_3^0 + \sum_{m=1}^2 \bar{u}_m^0 \dfrac{\partial h}{\partial \xi_m}
\right)\right]\nonumber\\
&&{}= -\dfrac{2}{\rho_0 \varepsilon h}  \bar{p}^2 + \nu  \left\{
\sum_{m=1}^2 \sum_{l=1}^2 \dfrac{
\partial^2   \bar{u}_3^1}{\partial \xi_l \partial \xi_m}
J^{0,0}_{lm} +  \sum_{k=1}^3 \sum_{l=1}^2\dfrac{
\partial  \bar{u}_k^1}{\partial \xi_l }
{{L}}^{1,0}_{kl3}(h)  +\sum_{k=1}^3 \bar{u}_k^1 S^{1,0}_{3k}(h)\right.\nonumber\\
&&\left.{}+\varepsilon \left[2 h \sum_{m=1}^2 \sum_{l=1}^2  \dfrac{
\partial^2  \bar{u}_3^0}{\partial \xi_l \partial \xi_m}
J^{1,0}_{lm}+ \sum_{k=1}^3 \sum_{l=1}^2\dfrac{
\partial \bar{u}_k^0}{\partial \xi_l } L_{kl3}^{0,1}(h) +  \sum_{k=1}^3 \bar{u}_k^0
S_{3k}^{0,1}(h)\right]
 \right.\nonumber\\
&&\left.{}+
\dfrac{2}{\varepsilon h}\dfrac{A^1}{A^0} \bar{u}_3^2 + \dfrac{6}{\varepsilon^2 h^2} \bar{u}_3^3\right\} + {\bar{f}_3^1}\label{ec_u31}
\end{eqnarray}
 where coefficients $B^1_{lk}$, $C^{1,0}_l$, $L^{1,0}_{kli}$, $L^{0,1}_{kli}$, $S^{1,0}_{ik}$, $S^{0,1}_{ik}$ are given by \eqref{B}, \eqref{Cij}, \eqref{L10}, \eqref{L01}-\eqref{L3_01}, \eqref{S10}-\eqref{S3_10}, \eqref{S01}-\eqref{S3_01} respectively. 

Analogously, we identify the terms multiplied by $\xi_3^2$ in \eqref{ec_ns_ij_alfa_beta_pol_xi3}-\eqref{div_i_dr_alfa_beta_pol_xi3} and, repeating the process once again, we obtain:
\begin{eqnarray}
&&\hspace*{-0.5cm}   \bar{u}_3^3=-\dfrac{\varepsilon h}{3}\left( \sum_{k=1}^2\dfrac{
\partial \bar{u}_k^2}{\partial \xi_k} + \sum_{k=1}^3 \sum_{l=1}^2 \bar{u}_k^2 H^0_{llk} - \dfrac{2}{h}  \sum_{k=1}^2  \bar{u}_k^2 \dfrac{\partial h}{\partial \xi_k} \right) \nonumber\\
&& {}-\dfrac{\varepsilon^2 h}{3} \left[ h \sum_{k=1}^2 \sum_{l=1}^2  \dfrac{
\partial \bar{u}_k^1}{\partial \xi_l} B^1_{lk} +\sum_{k=1}^2 \bar{u}_k^1 B^1_{3k}(h) +  h \sum_{k=1}^3 \bar{u}_k^1 \sum_{l=1}^2 H^1_{llk} \right]\nonumber\\
&& {}-\dfrac{\varepsilon^3 h^3}{3}\left[  \sum_{k=1}^2 \sum_{l=1}^2 \dfrac{
\partial \bar{u}_k^0}{\partial \xi_l}  B^2_{lk} +   \sum_{k=1}^3  \bar{u}_k^0 \sum_{l=1}^2 H^2_{llk} \right]
 \label{u33}\\
 &&\hspace*{-0.5cm}  \dfrac{ \partial \bar{u}_i^2}{\partial t} + \sum_{l=1}^2  \dfrac{
\partial \bar{u}_i^2}{\partial \xi_l} \left( \bar{u}_l^0 - C^0_l\right) +
\sum_{k=1}^3 \bar{u}_k^2 \left(Q^0_{ik}  + \sum_{l=1}^2 \bar{u}_l^0 H^0_{ilk}\right)\nonumber\\
&&{} + \dfrac{2}{h}  \bar{u}_i^2  \left[  \dfrac{1}{\varepsilon } \bar{u}_3^1  - \sum_{l=1}^2 \bar{u}_l^0 \dfrac{\partial h}{\partial \xi_l} -  \dfrac{\partial h}{\partial
t} -  C^0_3
\right]  +\dfrac{1}{\varepsilon h} \bar{u}_3^2  \bar{u}_i^1 \nonumber\\
&&{}+ \sum_{l=1}^2 \left(\dfrac{
\partial \bar{u}_i^0}{\partial \xi_l} +\sum_{k=1}^3 \bar{u}_k^0 H^0_{ilk} \right)  \left[ \bar{u}_l^2+ \varepsilon  h  \sum_{m=1}^2 \bar{u}_m^1
 B^1_{lm}  + (\varepsilon  h)^2 \left( \sum_{m=1}^2 u^0_m    B^2_{lm} -C^{2,1}_l  \right)\right]\nonumber\\
&&{}+ \sum_{l=1}^2 \left(   \dfrac{
\partial \bar{u}_i^1}{\partial \xi_l} + \sum_{k=1}^3 \bar{u}_k^1 H^0_{ilk}  \right) \left[ \bar{u}_l^1 -\varepsilon  h C^{1,0}_l +\varepsilon  h  \sum_{m=1}^2 u^0_m  B^1_{lm}  \right]\nonumber\\
&& {} + u^1_i \left( \varepsilon \sum_{l=1}^2 \bar{u}_l^0 B^1_{3l}  - \varepsilon C^{1,0}_3 -   \dfrac{1}{h}\sum_{l=1}^2 \bar{u}_l^1 \dfrac{\partial h}{\partial \xi_l} \right) + \dfrac{3}{\varepsilon h}   \bar{u}_i^3 \left[\bar{u}_3^0 - \left(\vec{a}_3 \cdot\dfrac{\partial
\vec{X}}{\partial t} \right)\right] \nonumber\\
&&=-\dfrac{1}{\rho_0} \left[ \sum_{l=1}^2\dfrac{
\partial  \bar{p}^2}{\partial \xi_l} J^{0,0}_{il}  +\varepsilon  h \sum_{l=1}^2\dfrac{
\partial  \bar{p}^1}{\partial \xi_l}  J^{0,1}_{il} + (\varepsilon  h)^2 \sum_{l=1}^2\dfrac{
\partial  \bar{p}^0}{\partial \xi_l} J^{0,2}_{il}   \right]  \nonumber\\
&&{}-\dfrac{2}{\rho_0 h}  \bar{p}^2  J^{0,0}_{i3} -\dfrac{ \varepsilon }{\rho_0}  \bar{p}^1 J^{0,1}_{i3}  + \nu \left\{ \dfrac{3}{\varepsilon h} 
\bar{u}_i^3\sum_{m=1}^2 H^0_{mm3}  +\sum_{l=1}^2 \sum_{m=1}^2  \dfrac{
\partial^2  \bar{u}_i^2}{\partial \xi_l \partial \xi_m} J^{0,0}_{lm} + \sum_{k=1}^3\sum_{l=1}^2 \dfrac{
\partial \bar{u}_k^2}{\partial \xi_l } {L}^{2,0}_{kli}(h)
\right.\nonumber\\
&&{}+   \sum_{k=1}^3 \bar{u}_k^2 {{S}^{2,0}_{ik}}(h)+   \varepsilon   \left[2 h  \sum_{l=1}^2 \sum_{m=1}^2 
\dfrac{
\partial^2   \bar{u}_i^1}{\partial \xi_l \partial \xi_m}  J^{1,0}_{lm}  +  \sum_{k=1}^3  \sum_{l=1}^2 \dfrac{
\partial \bar{u}_k^1}{\partial \xi_l } {L}^{1,1}_{kli}(h) 
+\sum_{k=1}^3 \bar{u}_k^1 S_{ik}^{1,1}(h)\right]\nonumber\\
&&\left.{}+ \varepsilon^2    \left[ \sum_{l=1}^2 \sum_{m=1}^2
\dfrac{
\partial^2   \bar{u}_i^0}{\partial \xi_l \partial \xi_m} \iota_{lm}^{2,1}(h) +    \sum_{l=1}^2  \sum_{k=1}^3 \dfrac{
\partial \bar{u}_k^0}{\partial \xi_l }{L}^{0,2}_{kli}(h)  + \sum_{k=1}^3 \bar{u}_k^0S_{ik}^{0,2}(h)\right]
\right\} + \bar{f}_i^2\quad (i=1,2) \label{ui_2}\\
&&\hspace*{-0.5cm} \dfrac{ \partial \bar{u}_3^2}{\partial t} + \sum_{l=1}^2  \dfrac{
\partial \bar{u}_3^2}{\partial \xi_l} \left( u^0_l-C^0_l  \right) +
\sum_{k=1}^2 \bar{u}_k^2\left[ Q^0_{3k}+ \sum_{l=1}^2 \bar{u}_l^0 \left( \vec{a}_3 \cdot\dfrac{
\partial \vec{a}_{k}}{\partial \xi_l}\right)\right]\nonumber\\
&&{} +\dfrac{3}{\varepsilon h} \bar{u}_3^2 \bar{u}_3^1 - \dfrac{2}{h} \bar{u}_3^2\left[ \sum_{l=1}^2 \bar{u}_l^0 \dfrac{\partial h}{\partial \xi_l}+ \dfrac{\partial h}{\partial
t}+  C^0_3
\right]\nonumber\\
&&{}+ \sum_{l=1}^2 \left( \dfrac{
\partial \bar{u}_3^0}{\partial \xi_l}+ \sum_{k=1}^2 \bar{u}_k^0  \left( \vec{a}_3 \cdot\dfrac{
\partial \vec{a}_{k}}{\partial \xi_l}\right)  \right)  \left[ \bar{u}_l^2+ \varepsilon  h  \sum_{m=1}^2 \bar{u}_m^1
B^1_{lm}  - (\varepsilon  h)^2 C^{2,1}_l + (\varepsilon h)^2 \sum_{m=1}^2 \bar{u}_m^0   B^2_{lm}\right]\nonumber\\
&&{}+ \sum_{l=1}^2 \left(  \dfrac{
\partial \bar{u}_3^1}{\partial \xi_l} + \sum_{k=1}^2 \bar{u}_k^1  \left( \vec{a}_3 \cdot\dfrac{
\partial \vec{a}_{k}}{\partial \xi_l}\right)  \right) \left[ \bar{u}_l^1-\varepsilon  h C^{1,0}_l  + \varepsilon  h   \sum_{m=1}^2 \bar{u}_m^0   B^1_{lm}\right]    \nonumber\\
&&{} +  \bar{u}_3^1\left(\varepsilon \sum_{l=1}^2 \bar{u}_l^0  B^1_{3l} - \varepsilon C^{1,0}_3   - \dfrac{1}{h} \sum_{l=1}^2 \bar{u}_l^1 \dfrac{\partial h}{\partial \xi_l}\right)  + \dfrac{3}{\varepsilon h}   \bar{u}_3^3 \left[\bar{u}_3^0 - \left(\vec{a}_3 \cdot\dfrac{\partial
\vec{X}}{\partial t} \right)\right]  \nonumber\\
&&=-\dfrac{3}{\varepsilon \rho_0 h} \bar{p}^3 + \nu \left\{\dfrac{3}{\varepsilon h} \bar{u}_3^3\sum_{m=1}^2 H^0_{mm3}+ \sum_{l=1}^2 \sum_{m=1}^2  \dfrac{
\partial^2  \bar{u}_3^2}{\partial \xi_l \partial \xi_m} J^{0,0}_{lm}  + \sum_{k=1}^3 \sum_{l=1}^2 
 \dfrac{
\partial \bar{u}_k^2}{\partial \xi_l } L^{2,0}_{kl3}(h) + \sum_{k=1}^3 \bar{u}_k^2S^{2,0}_{3k}(h)\right.
\nonumber\\
&&{}+ \varepsilon \left[ 2 h \sum_{l=1}^2 \sum_{m=1}^2 \dfrac{
\partial^2   \bar{u}_3^1}{\partial \xi_l \partial \xi_m}  J^{1,0}_{lm} + \sum_{l=1}^2  \sum_{k=1}^3\dfrac{
\partial \bar{u}_k^1}{\partial \xi_l } L_{kl3}^{1,1}(h)  +
\sum_{k=1}^3 \bar{u}_k^1 S_{3k}^{1,1}(h)\right]\nonumber\\
&&\left.{}+ \varepsilon^2   \left[ \sum_{l=1}^2 \sum_{m=1}^2 
\dfrac{
\partial^2   \bar{u}_3^0}{\partial \xi_l \partial \xi_m} \iota^{2,1}_{lm}  +  \sum_{k=1}^3 \sum_{l=1}^2  \dfrac{
\partial \bar{u}_k^0}{\partial \xi_l } {L}^{0,2}_{kl3}(h)+ \sum_{k=1}^3 \bar{u}_k^0 S_{3k}^{0,2}(h)\right] 
\right\} + \bar{f}_3^2\label{u3_2}
\end{eqnarray}
where the coefficients that appear are defined in appendix \ref{ApendiceA} (see \eqref{B}-\eqref{Cij}, \eqref{H}, \eqref{J}, \eqref{Q}, \eqref{Q3}, \eqref{L20}-\eqref{L3_02}, \eqref{S20}-\eqref{S3_02}, \eqref{iota21}).

And, finally, we identify the terms multiplied by $\xi_3^3$  in \eqref{ec_ns_ij_alfa_beta_pol_xi3} - \eqref{div_i_dr_alfa_beta_pol_xi3} and, following the same steps we have:
\begin{eqnarray}
&&\hspace*{-0.5cm} \sum_{k=1}^2\dfrac{
\partial \bar{u}_k^3}{\partial \xi_k}  + \sum_{k=1}^3  \bar{u}_k^3 \sum_{l=1}^2  H^0_{llk} -\dfrac{ 3}{h} \sum_{k=1}^2  \bar{u}_k^3 \dfrac{\partial h}{\partial \xi_k} \nonumber\\
&&{} +  \varepsilon \left[  h\sum_{k=1}^2  \sum_{l=1}^2 \dfrac{
\partial \bar{u}_k^2}{\partial \xi_l}   B^1_{lk} + h \sum_{k=1}^3 \bar{u}_k^2 \sum_{l=1}^2 H^1_{llk}  + 2  \sum_{k=1}^2 \bar{u}_k^2 B^1_{3k}(h) \right]\nonumber\\
&&{}+ \varepsilon^2\left[  h^2  \sum_{k=1}^2  \sum_{l=1}^2 \dfrac{
\partial \bar{u}_k^1}{\partial \xi_l}   B^2_{lk} + h^2 \sum_{k=1}^3 \bar{u}_k^1 \sum_{l=1}^2 H^2_{llk}+ 
 h  \sum_{k=1}^2 \bar{u}_k^1 B^2_{3k}(h)\right]\nonumber\\
&&{} +  (\varepsilon h)^3  \left\{\sum_{k=1}^2  \sum_{l=1}^2 \dfrac{
\partial \bar{u}_k^0}{\partial \xi_l} B^3_{lk}+  \sum_{k=1}^3 \bar{u}_k^0 \sum_{l=1}^2 H^3_{llk}\right\}  =0\nonumber\\
\label{div_i_dr_alfa_beta_pol_xi3_3s}
\\\
&&\hspace*{-0.5cm} \dfrac{ \partial \bar{u}_i^3}{\partial t} +   \sum_{l=1}^2  \dfrac{
\partial \bar{u}_i^3}{\partial \xi_l} (u^0_l -C^0_l)+ \sum_{k=1}^3
\bar{u}_k^3 \left[ Q^0_{ik} +  \sum_{l=1}^2 H^0_{ilk} u^0_l \right]   +\sum_{k=1}^2 \bar{u}_k^3\left(  \dfrac{
\partial \bar{u}_i^0}{\partial \xi_k} +  \sum_{m=1}^3 \bar{u}_m^0 H^0_{ikm}\right)\nonumber\\
&&{}+\dfrac{3}{ h} \bar{u}_i^3\left( \dfrac{1}{\varepsilon}\bar{u}_3^1 -\dfrac{\partial h}{\partial
t}- C_3^0 - \sum_{m=1}^2 \bar{u}_m^0 \dfrac{\partial h}{\partial \xi_m}
\right)  +\dfrac{1}{\varepsilon h}  \bar{u}_3^3  \bar{u}_i^1 \nonumber\\
&&{}+ \sum_{l=1}^2   \dfrac{
\partial \bar{u}_i^2}{\partial \xi_l} \left(\bar{u}_l^1- \varepsilon hC^{1,0}_l  + \varepsilon h  \sum_{m=1}^2 \bar{u}_m^0 B^1_{lm}\right)  \nonumber\\
&&{}+ \sum_{k=1}^3 \bar{u}_k^2 \sum_{l=1}^2 H^0_{ilk} \left(u^1_l - \varepsilon h C^{1,0}_l+ \varepsilon h \sum_{m=1}^2 \bar{u}_m^0  B^1_{lm} \right)\nonumber\\
&&{}+\sum_{k=1}^2 \bar{u}_k^2 \left[ \dfrac{
\partial \bar{u}_i^1}{\partial \xi_k} -\dfrac{1}{h} \dfrac{\partial h}{\partial \xi_k}  \bar{u}_i^1 +\sum_{m=1}^3   \bar{u}_m^1H^0_{ikm} + \varepsilon h \sum_{l=1}^2\left(  \dfrac{
\partial \bar{u}_i^0}{\partial \xi_l} + \sum_{m=1}^3  \bar{u}_m^0 H^0_{ikm} \right)  B^1_{lm} \right]\nonumber\\
&&{}+ 2 \bar{u}_i^2 \left( \dfrac{1}{\varepsilon h} \bar{u}_3^2     -\dfrac{1}{h} \sum_{k=1}^2 \bar{u}_k^1\dfrac{\partial h}{\partial \xi_k}  - \varepsilon C^{1,0}_3 +  \varepsilon \sum_{k=1}^2 \bar{u}_k^0 
B^1_{3k}(h)  \right) \nonumber\\
&&{}+\varepsilon h  \sum_{l=1}^2 \left( \dfrac{
\partial \bar{u}_i^1}{\partial \xi_l} + \sum_{k=1}^3 \bar{u}_k^1H^0_{ilk}\right)  \left[ \sum_{m=1}^2 \left( \bar{u}_m^1  B^1_{lm} + \varepsilon  h \bar{u}_m^0 B^2_{lm}\right) - \varepsilon  h  C^{2,1}_l  \right]\nonumber\\
&&{} +\varepsilon \bar{u}_i^1\sum_{k=1}^2 \left[ \bar{u}_k^1 B^1_{3k}(h) +  \varepsilon h     \bar{u}_k^0
B^2_{3k}(h) - \varepsilon h     C^{2,1}_3 \right]\nonumber\\
&&{} + (\varepsilon  h)^2 \sum_{k=1}^2 \bar{u}_k^1  \sum_{l=1}^2\left(  \dfrac{
\partial \bar{u}_i^0}{\partial \xi_l} + \sum_{m=1}^3  \bar{u}_m^0 H^0_{ilm}\right) B^2_{lk}\nonumber\\
&&{}+ (\varepsilon  h)^3 \sum_{k=1}^3 \sum_{l=1}^2 \left( \dfrac{
\partial \bar{u}_i^0}{\partial \xi_l} +  \bar{u}_k^0H^0_{ilk} \right) \left[\sum_{m=1}^2 \bar{u}_m^0 B^3_{lm} -C^{3,2}_l \right]
\nonumber\\
&&=-\dfrac{1}{\rho_0} \sum_{l=1}^2 \left[  \dfrac{
\partial  \bar{p}^3}{\partial \xi_l} J^{0,0}_{il} +\varepsilon h  \dfrac{
\partial  \bar{p}^2}{\partial \xi_l} J^{0,1}_{il}  + (\varepsilon  h)^2 \dfrac{
\partial  \bar{p}^1}{\partial \xi_l}  J^{0,2}_{il} + (\varepsilon h)^3 \dfrac{
\partial  \bar{p}^0}{\partial \xi_l} J^{0,3}_{il} \right]\nonumber\\
&&{} -\dfrac{1}{\rho_0} \left[\dfrac{3}{h} \bar{p}^3 J^{0,0}_{i3}+ 2 \varepsilon  \bar{p}^2 J^{0,1}_{i3} + \varepsilon^2 h \bar{p}^1 J^{0,2}_{i3}\right]      \nonumber\\
&&{} + \nu \left[\sum_{m=1}^2 \sum_{l=1}^2 \dfrac{
\partial^2  \bar{u}_i^3}{\partial \xi_l \partial \xi_m}
J^{0,0}_{lm} + \sum_{k=1}^3 \sum_{l=1}^2 \dfrac{
\partial  \bar{u}_k^3}{\partial \xi_l} {L}^{3,0}_{kli}(h) + \sum _{k=1}^3 \bar{u}_k^3 {{S}^{3,0}_{ik}}(h) \right. \nonumber\\
&&\left.{}+2 \varepsilon h  \sum_{m=1}^2 \sum_{l=1}^2 \dfrac{
\partial^2  \bar{u}_i^2}{\partial \xi_l \partial \xi_m}
J^{1,0}_{lm} + \varepsilon   \sum_{k=1}^3  \sum_{l=1}^2   \dfrac{
\partial   \bar{u}_k^2}{\partial \xi_l} {L}^{2,1}_{kli}(h) +\varepsilon \sum_{k=1}^3 \bar{u}_k^2S_{ik}^{2,1}(h) \right. \nonumber \\
&&\left.{}+ \varepsilon^2   \sum_{l=1}^2 \sum_{m=1}^2\dfrac{
\partial^2   \bar{u}_i^1}{\partial \xi_l \partial \xi_m} \iota^{2,1}_{lm}(h) +\varepsilon^2 \sum_{k=1}^3
  \sum_{l=1}^2 \dfrac{
\partial  \bar{u}_k^1}{\partial \xi_l}{L}^{1,2}_{kli}(h) + \varepsilon^2 \sum_{k=1}^3  \bar{u}_k^1 S_{ik}^{1,2}(h) \right. \nonumber\\
&&\left.{} +  \varepsilon^3 \sum_{m=1}^2 \sum_{l=1}^2\dfrac{
\partial^2   \bar{u}_1^0}{\partial \xi_l \partial \xi_m}
\iota^3_{lm}(h) +  \varepsilon^3  \sum_{k=1}^3 \sum_{l=1}^2\dfrac{
\partial  \bar{u}_k^0}{\partial \xi_l }{L}^{0,3}_{kli}(h)+  \varepsilon^3   \sum_{k=1}^3 \bar{u}_k^0 S_{ik}^{0,3}(h)
\right]+
\bar{f}_i^3, \quad (i=1,2) \nonumber \\
&& \label{ui_3}
\end{eqnarray} 
where the coefficients that appear are defined in appendix \ref{ApendiceA} (see \eqref{B}-\eqref{Cij}, \eqref{H}, \eqref{J}, \eqref{Q}, \eqref{L30}, \eqref{L21}, \eqref{L12}, \eqref{L03}, \eqref{S30}, \eqref{S21}, \eqref{S12}, \eqref{S03}, \eqref{iota21}, \eqref{iota3}).

An equation analogous to \eqref{ec_ns_u30_Res}, \eqref{ec_u31} and \eqref{u3_2} is also obtained when identifying the terms multiplied by $\xi_3^3$ in \eqref{ec_ns_ij_alfa_beta_pol_xi3} and multiplying by $\vec{a}_3$. We do not write it explicitly below because, as we will see later, equation \eqref{u30} makes it unnecessary.

Since we have assumed that the velocity and the pressure are polynomials of degree three in $\xi_3$ (\eqref{ui_pol_xi3}-\eqref{p_pol_xi3}), we have 16 unknowns to determine. Out of these unknowns,  the terms $\bar{u}^k_3$ and $\bar{p}^k$ ($k=1,2,3$) corresponding to the third component of the velocity and the pressure, respectively, are given by \eqref{u_3^1}, \eqref{u32}, \eqref{u33}, \eqref{ec_ns_u30_Res}, \eqref{ec_u31} and \eqref{u3_2} using the terms $\bar{u}^k_i$ ($i=1,2, k=0,1,2$) once they have been computed. Therefore, we must actually determine 10 unknowns.

Now, if we denote by $V_1 \vec{a}_1+V_2 \vec{a}_2$ the tangential velocity
at the lower surface and by $W_1\vec{a}_1+W_2\vec{a}_2$ the
tangential velocity at the upper surface, we have
\begin{eqnarray}
u_k^{\varepsilon} \vec{e}_k = u_k(\varepsilon)\vec{a}_{k}&=&
V_1(\varepsilon) \vec{a}_1+V_2(\varepsilon) \vec{a}_2+\left( \dfrac{\partial
\vec{X}}{\partial t}
\cdot \vec{a}_3\right)\vec{a}_3 \ \textrm{on }\xi_3=0 \label{cc_xi3_0}\\
u_k^{\varepsilon} \vec{e}_k =u_k(\varepsilon)\vec{a}_{k}&=&
W_1(\varepsilon) \vec{a}_1+ W_2(\varepsilon) \vec{a}_2+\left( \dfrac{\partial (\vec{X}+
\varepsilon h\vec{a}_3)}{\partial t} \cdot \vec{a}_3\right)\vec{a}_3
\ \textrm{on }\xi_3=1 \label{cc_xi3_1}
\end{eqnarray}
and, taking into account \eqref{ui_pol_xi3}, we yield
\begin{eqnarray}
&&\bar{u}_i^0=V_i \quad (i=1,2) \label{ui0_Vi}\\
&&\bar{u}_3^0= \dfrac{\partial \vec{X}}{\partial t}
\cdot \vec{a}_3 \label{u30}\\
&&\sum_{k=1}^3  \bar{u}_i^k= W_i - V_i \quad (i=1,2)
 \label{cc_uk_xi3_1b}\\
&&\sum_{k=1}^3  \bar{u}_3^k=  \varepsilon \dfrac{\partial
h}{\partial t}
\label{cc_u3_xi3_1b}
\end{eqnarray} 

Equality \eqref{u30} gives us an expression for $\bar{u}_3^0$, so it is no longer an unknown, it is determined by the lower bound surface. At this point, 9 unknowns are left, $\bar{u}_i^k$ $(i=1,2,\ k=0,1,2,3)$ and $\bar{p}^0$, but we will see that not all are needed to obtain an approximation of the velocity and the pressure.

\subsection{Stating the order in $\varepsilon$ of the unknowns}

We can assume, after \eqref{ansatz_1}, that the velocity is of order $\varepsilon^0$, but from \eqref{u_3^1},  \eqref{u32} and \eqref{u33} we know that, actually, the terms $\bar{u}_3^k$, ($k=1,2,3$) are of order $\varepsilon$.

Regarding the pressure terms, we deduce from equations \eqref{ec_ui0_coefs} that $\bar{p}^0$ is of order $\varepsilon^{-2}$, which fits with  hypothesis \eqref{ansatz_2}. 
The term  $\bar{p}^1$ can be obtained from \eqref{ec_ns_u30_Res} that, using \eqref{u30}, \eqref{u_3^1}, \eqref{u32}, writes
\begin{eqnarray}
&&\hspace*{-0.5cm} \dfrac{1}{\rho_0} \bar{p}^1=-\varepsilon h\dfrac{ \partial \bar{u}_3^0}{\partial t} -\varepsilon h\sum_{l=1}^2  \dfrac{ \partial \bar{u}_3^0}{\partial \xi_l} (\bar{u}_l^0 -C^0_l)
 -\varepsilon h \sum_{k=1}^2
 \bar{u}_k^0\left[ Q^0_{3k}+ \sum_{l=1}^2 \bar{u}_l^0 \left(\vec{a}_3 \cdot \dfrac{ \partial \vec{a}_{l}}{\partial \xi_k}\right) \right]
\nonumber\\
&&{}  + \nu \left[ \varepsilon h\sum_{m=1}^2 \sum_{l=1}^2\dfrac{
\partial^2   \bar{u}_3^0}{\partial \xi_l \partial \xi_m}
J^{0,0}_{lm}  +\varepsilon h\sum_{k=1}^3 \sum_{l=1}^2 \dfrac{
\partial  \bar{u}_k^0}{\partial \xi_l } L^0_{kl3} + \varepsilon h\sum_{k=1}^3  \bar{u}_k^0 {{S}^0_{3k}}\right. \nonumber\\
&&\left.{} - \dfrac{1}{\sqrt{A^0}} \textrm{div}(\sqrt{A^0}\vec{u}^1) + \dfrac{1 }{h} \nabla h \cdot \vec{u}^1 + \varepsilon h\left( \dfrac{1}{\sqrt{A^0}} \textrm{div}(\sqrt{A^0}\vec{u}^0)  {+ \bar{u}_3^0 \dfrac{A^1}{A^0}}\right)\dfrac{A^1}{A^0}  \right.\nonumber\\
&&\left.{} -\varepsilon h \sum_{l=1}^2  \sum_{k=1}^3 \left[ \dfrac{
\partial \bar{u}_k^0}{\partial \xi_l}  B^1_{lk}+   \bar{u}_k^0  H^1_{llk}\right]
\right]
+\varepsilon h {\bar{f}_3^0}\label{p1}
\end{eqnarray} 

So we have,
\begin{eqnarray}
&&\hspace*{-0.5cm}  \bar{p}^1=\mu \left( - \dfrac{1}{\sqrt{A^0}} \textrm{div}(\sqrt{A^0}\vec{u}^1) + \dfrac{1 }{h} \nabla h \cdot \vec{u}^1 \right) + O(\varepsilon)\label{p10}
\end{eqnarray} 

And finally, for the terms $\bar{p}^2$ and $\bar{p}^3$, from \eqref{ec_u31} and \eqref{u3_2}, using \eqref{u30}, \eqref{u_3^1}, \eqref{u32}, \eqref{u33}, we get
\begin{eqnarray}
&&\hspace*{-0.5cm}  \bar{p}^2 =- \mu \left( \sum_{k=1}^2\dfrac{
\partial \bar{u}_k^2}{\partial \xi_k} + \sum_{k=1}^3 \sum_{l=1}^2 \bar{u}_k^2 H^0_{llk} - \dfrac{2}{h}  \sum_{k=1}^2  \bar{u}_k^2 \dfrac{\partial h}{\partial \xi_k} \right) +O(\varepsilon)\label{p20}
\\
&&\hspace*{-0.5cm}  \bar{p}^3= -\dfrac{\varepsilon \rho_0 h}{3}\left\{
\sum_{k=1}^2 \bar{u}_k^2\left[ Q^0_{3k} + \dfrac{
\partial \bar{u}_3^0}{\partial \xi_k}+ \sum_{l=1}^2 \bar{u}_l^0 \left( \vec{a}_3 \cdot\dfrac{
\partial \vec{a}_{k}}{\partial \xi_l}+\vec{a}_3 \cdot\dfrac{
\partial \vec{a}_{l}}{\partial \xi_k}\right)  \right] +  \sum_{l=1}^2  \bar{u}_l^1 \sum_{k=1}^2 \bar{u}_k^1  \left( \vec{a}_3 \cdot\dfrac{
\partial \vec{a}_{k}}{\partial \xi_l}\right)  \right\}  \nonumber\\
&&{} + \dfrac{\varepsilon \mu h}{3} \left\{
-\left( \sum_{k=1}^2\dfrac{
\partial \bar{u}_k^2}{\partial \xi_k} + \sum_{k=1}^3 \sum_{l=1}^2 \bar{u}_k^2 H^0_{llk} - \dfrac{2}{h}  \sum_{k=1}^2  \bar{u}_k^2 \dfrac{\partial h}{\partial \xi_k} \right) 
\sum_{m=1}^2 H^0_{mm3} \right.\nonumber\\
&&\left.{}+\sum_{k=1}^2 \sum_{l=1}^2 
 \dfrac{
\partial \bar{u}_k^2}{\partial \xi_l } L^{2,0}_{kl3}(h) +\sum_{k=1}^2 \bar{u}_k^2S^{2,0}_{3k}(h)
\right\} +\dfrac{\varepsilon \rho_0 h}{3} \bar{f}_3^2 +O(\varepsilon^2) \label{p3e1}
\end{eqnarray}

We also obtain that the terms $\bar{u}_1^3$ and $\bar{u}_2^3$ are of order $\varepsilon$, since equations \eqref{ec_u1i} could be written, taking into account \eqref{u_3^1}, \eqref{p10} and that $\bar{p}^0$ is of order $\varepsilon^{-2}$, as follows
\begin{eqnarray}
&&\hspace*{-0.5cm}
\bar{u}_i^3 =-  \dfrac{\varepsilon h}{ 3} \dfrac{A^1}{A^0} \bar{u}_i^2  - \dfrac{\varepsilon^3 h}{\rho_0} \sum_{l=1}^2  \dfrac{
\partial  \bar{p}^0}{\partial \xi_l} J^{0,1}_{il} + O(\varepsilon^2) \quad (i=1,2)\label{ui3_ui2}
\end{eqnarray}

\section{Imposing boundary conditions} \label{BoundaryConditions}
\subsection{Boundary conditions leading to a lubrication problem} \label{subseccion-4-1}

Let us assume that the fluid slips at the lower
surface $(\xi _{3}=0)$, and at the upper surface $(\xi _{3}=1)$, that is, 
let us assume that the tangential velocities of the fluid at
the lower and upper surfaces are known, and, that the normal velocity of both surfaces matches the normal velocities of the fluid at the surfaces.

In this case, the terms $\bar{u}_i^0$ ($i=1,2$) are known (\eqref{ui0_Vi}), and from equations \eqref{ec_ui0_coefs} we obtain
\begin{equation}
\bar{u}_i^{2}=\dfrac{h^2 \varepsilon^2}{2\mu}\sum_{l=1}^2 \dfrac{\partial \bar{p}^{0}}{\partial \xi_l} J^{0,0}_{il}+O(\varepsilon) \label{u_i^2_p0},\quad  (i=1,2)
\end{equation}

Then, considering \eqref{ui3_ui2}, \eqref{cc_uk_xi3_1b} writes: 
\begin{equation}
    \bar{u}_i^1=W_i-V_i-\bar{u}_i^2+O(\varepsilon) \label{ui1_WVui2}
\end{equation}

We can substitute  $\bar{u}^k_3$ in \eqref{cc_u3_xi3_1b}  by the expressions \eqref{u_3^1}, \eqref{u32} and \eqref{u33}:
\begin{eqnarray}
&&\hspace{-0.5cm} {} -\dfrac{ h}{\sqrt{A^0}} \textrm{div}(\sqrt{A^0}\vec{V})
- h \left( \dfrac{\partial \vec{X}}{\partial t}
\cdot \vec{a}_3 \right)\dfrac{A^1}{A^0}- \dfrac{ h}{2\sqrt{A^0}} \textrm{div}(\sqrt{A^0}\vec{u}^1) + \dfrac{1 }{2} \nabla h \cdot \vec{u}^1 \nonumber\\
&&{}-\dfrac{ h}{3}\sum_{k=1}^2\dfrac{
\partial \bar{u}_k^2}{\partial \xi_k} -\dfrac{ h}{3}\sum_{k=1}^3 \bar{u}_k^2\sum_{l=1}^2  H^0_{llk} + \dfrac{2  }{3} \sum_{k=1}^2  \bar{u}_k^2 \dfrac{\partial h}{\partial \xi_k} +O(\varepsilon)=   \dfrac{\partial
h}{\partial t}\label{ec_previa_lub}
\end{eqnarray}
and, next, rewrite, using expressions \eqref{u_i^2_p0} and \eqref{ui1_WVui2}:

\begin{eqnarray}
&&\hspace*{-0.5cm} {}- h \left( \dfrac{\partial \vec{X}}{\partial t}
\cdot \vec{a}_3 \right)\dfrac{A^1}{A^0} - \dfrac{ h}{2\sqrt{A^0}} \textrm{div}(\sqrt{A^0}(\vec{W}+\vec{V})) + \dfrac{1 }{2} \nabla h \cdot (\vec{W}-\vec{V})\nonumber\\
&&{}+ \dfrac{ \varepsilon^2}{12 \mu\sqrt{A^0}} \textrm{div} \left(h^3\sqrt{A^0}\sum_{l=1}^2 \dfrac{\partial \bar{p}^{0}}{\partial \xi_l} \left(J^{0,0}_{1l}, J^{0,0}_{2l}\right) \right) =   \dfrac{\partial
h}{\partial t} +O(\varepsilon)\label{lubric}
\end{eqnarray}

Taking into account that 
$$ \dfrac{M}{\sqrt{A^0}}= \sqrt{A^0}  \begin{pmatrix}
J_{11}^{0,0}& J_{12}^{0,0}\\ J_{21}^{0,0} & J_{22}^{0,0}
 \end{pmatrix} $$
 and that $\bar{p}^{0}$ is of order $\varepsilon^{-2}$, we obtain that equation \eqref{lubric} is equivalent to \eqref{Reynolds_gen}.

Once we solve equation  \eqref{lubric} to obtain $\bar{p}^0$, using \eqref{ui0_Vi}, \eqref{ui1_WVui2}, \eqref{u_i^2_p0}, \eqref{ui3_ui2}, \eqref{u30}, \eqref{u_3^1}, \eqref{u32}, \eqref{u33}, \eqref{ec_ui0_coefs}, \eqref{p10}, \eqref{p20} and \eqref{p3e1}, we have 
\begin{eqnarray}
u_i^{\varepsilon}&=& \dfrac{(h^\varepsilon)^2 (\xi_3^2-\xi_3)}{2\mu}\sum_{l=1}^2 \dfrac{\partial \bar{p}^{0}}{\partial \xi_l} J^{0,0}_{il}+\left( W_i-V_i\right)\xi_3+V_i+O(\varepsilon), \quad (i=1,2) \label{eq-u-i-4-1} \\
u_3^{\varepsilon}&=&  \dfrac{\partial \vec{X}}{\partial t}
\cdot \vec{a}_3 +O(\varepsilon)\\
{p}^{\varepsilon}&=& \bar{p}^0 +\xi_3 \bar{p}^1+\xi_3^2 \bar{p}^2+O(\varepsilon)
\end{eqnarray}

\subsection{Boundary conditions leading to a shallow water problem} \label{subseccion-4-2}
 
Now, instead
of considering that the tangential and normal velocities are known on the
upper and lower surfaces,
we assume that the normal component of the traction on $\xi _{3}=0$ and
on $\xi _{3}=1$ are known pressures, and that the tangential component
of the traction on these surfaces are friction forces depending on the
value of the velocities on $\partial D$. Therefore, we assume that
\begin{eqnarray}
&&\vec{T}^{\varepsilon}\cdot
\vec{n}^{\varepsilon}_0 = (\sigma^{\varepsilon} \vec{n}^{\varepsilon}_0)\cdot
\vec{n}^{\varepsilon}_0=-\pi^{\varepsilon}_0 \textrm{ on } \xi_3=0, \label{Tn0_xi3_0}\\
&&\vec{T}^{\varepsilon}\cdot
\vec{n}^{\varepsilon}_1 = (\sigma^{\varepsilon} \vec{n}^{\varepsilon}_1)\cdot
\vec{n}^{\varepsilon}_1=-\pi^{\varepsilon}_1 \textrm{ on } \xi_3=1,
\label{Tn1_xi3_1}
\\
&&\vec{T}^{\varepsilon}\cdot
\vec{a}_i= (\sigma^{\varepsilon}  \vec{n}^{\varepsilon}_0)\cdot
\vec{a}_i=-\vec{f}^{\varepsilon}_{R_0}\cdot
\vec{a}_i \textrm{ on } \xi_3=0,\quad (i=1,2) \label{Tai_xi3_0}\\
&&\vec{T}^{\varepsilon}\cdot
\vec{v}_i^{\varepsilon}= (\sigma^{\varepsilon}   \vec{n}^{\varepsilon}_1)\cdot
\vec{v}_i^{\varepsilon}=-\vec{f}^{\varepsilon}_{R_1}\cdot
\vec{v}_i^{\varepsilon} \textrm{ on } \xi_3=1,\quad (i=1,2)
\label{Tvi_xi3_1}
\end{eqnarray}where $\vec{T}^{\varepsilon}$ is the traction vector and $\sigma^{\varepsilon}$ is the stress tensor given by
\begin{eqnarray}
    \sigma^{\varepsilon}_{ij}&=& -p^{\varepsilon}\delta_{ij}+\mu \left(
\dfrac{\partial u_i^{\varepsilon}}{\partial x^{\varepsilon}_j} +
\dfrac{\partial u_j^{\varepsilon}}{\partial
x^{\varepsilon}_i}\right)\nonumber\\
&=&   \sum_{k=0}^3 \xi_3^k \left[ -\bar{p}^k \delta_{ij}+\mu \sum_{m=1}^3 \sum_{l=1}^2 \left(
\dfrac{\partial( \bar{u}_m^k  a_{mi})}{\partial \xi_l} \left(\sum_{n=0}^{\infty}(\varepsilon \xi_3 h)^n(\alpha_l^n a_{1j}+\beta_l^n a_{2j})\right) \right. \right.\nonumber\\
&+&\left. \left.
\dfrac{\partial (\bar{u}_m^k a_{mj})}{\partial
\xi_l}  \left(\sum_{n=0}^{\infty}(\varepsilon \xi_3 h)^n(\alpha_l^n a_{1i}+\beta_l^n a_{2i})\right) \right)\right]\nonumber\\
&+& \mu  \sum_{k=1}^3 \sum_{m=1}^3 k \xi_3^{k-1}  \bar{u}_m^k \left[\sum_{n=0}^{\infty}\varepsilon^n \xi_3^{n+1} h^{n-1}  
 \left(\alpha_3^n(  a_{mi} a_{1j}+ a_{mj} a_{1i})+\beta_3^n (a_{mi} a_{2j}+ a_{mj} a_{2i}) \right) \right.\nonumber\\
&+&\left.\dfrac{1}{\varepsilon h} ( a_{mi}  a_{3j} +  a_{mj}a_{3i})  \right],  \quad (i,j=1,2,3) \label{sigmaij_des}
\end{eqnarray}
vectors $\vec{n}^{\varepsilon}_0$, $\vec{n}^{\varepsilon}_1$ are, respectively, the outward unit normal vectors to the lower and the upper surfaces, that is 
\begin{eqnarray}
&&\vec{n}^{\varepsilon}_0=s_0 \vec{a}_3 \label{n0}\\
&&\vec{n}^{\varepsilon}_1=-s_0\dfrac{\vec{v}^{\varepsilon}_3}{\|\vec{v}^{\varepsilon}_3\|}
\label{n1} 
\end{eqnarray}
where
\begin{equation}
    s_0=-1 \textrm{ or } s_0=1 
\end{equation}
is fixed ($\vec{n}^{\varepsilon}_0 = \vec{a}_3$ or $\vec{n}^{\varepsilon}_0 = - \vec{a}_3$, depending on the orientation of the parametrization $\vec{X}$), and 
\begin{eqnarray}
&&\vec{v}^{\varepsilon}_1=\vec{a}_1+\varepsilon
\left(\dfrac{\partial h}{\partial \xi_1} \vec{a}_3 + h
\dfrac{\partial \vec{a}_3}{\partial \xi_1}\right) \label{v1}\\
&&\vec{v}^{\varepsilon}_2=\vec{a}_2+\varepsilon
\left(\dfrac{\partial h}{\partial \xi_2} \vec{a}_3 + h
\dfrac{\partial \vec{a}_3}{\partial \xi_2}\right) \label{v2}\\
&&\vec{v}^{\varepsilon}_3=\vec{v}^{\varepsilon}_1 \times \vec{v}^{\varepsilon}_2 =\vec{a}_1 \times \vec{a}_2 +\varepsilon
\left[\dfrac{\partial h}{\partial \xi_2} (\vec{a}_1 \times
\vec{a}_3) + h \left(\vec{a}_1 \times \dfrac{\partial
\vec{a}_3}{\partial \xi_2}\right) + \dfrac{\partial h}{\partial
\xi_1} (\vec{a}_3\times \vec{a}_2) + h \left(\dfrac{\partial
\vec{a}_3}{\partial \xi_1} \times
\vec{a}_2\right)\right]\nonumber\\
&&\hspace*{0.5cm}+\varepsilon^2 \left[ \left(\dfrac{\partial h}{\partial \xi_1}
\vec{a}_3 + h \dfrac{\partial \vec{a}_3}{\partial
\xi_1}\right)\times\left(\dfrac{\partial h}{\partial \xi_2}
\vec{a}_3 + h \dfrac{\partial \vec{a}_3}{\partial \xi_2}\right)
\right]\label{v3}\\
 &&\|\vec{v}^{\varepsilon}_3\|=\|\vec{a}_1 \times
\vec{a}_2\| + \varepsilon h  \left[ \vec{a}_3 \cdot \left(\vec{a}_1
\times \dfrac{\partial \vec{a}_3}{\partial \xi_2}\right) +\vec{a}_3
\cdot \left(\dfrac{\partial \vec{a}_3}{\partial \xi_1} \times
\vec{a}_2\right)\right] + O(\varepsilon^2)\label{mod_v3}
\end{eqnarray}

Now, condition \eqref{Tn0_xi3_0} can be written (using
\eqref{sigmaij_des}, \eqref{n0}) as:
\begin{eqnarray}
&&(\sigma_{ij}a_{3j})a_{3i} = -\bar{p}^0 + \dfrac{2\mu }{\varepsilon h}  \bar{u}_3^1=-\pi_0
\end{eqnarray}
where $\pi_0 = \pi_0(\varepsilon) = \pi_0^\varepsilon$.

Taking into account \eqref{u_3^1} and \eqref{ui0_Vi}, we have
\begin{eqnarray}
\bar{p}^0&=& \dfrac{2\mu \bar{u}_3^1}{\varepsilon h} +\pi_0= -\dfrac{2\mu}{\sqrt{A^0}} \textrm{div}(\sqrt{A^0}\vec{V}){-2\mu\left( \dfrac{\partial \vec{X}}{\partial t}
\cdot \vec{a}_3 \right)\dfrac{A^1}{A^0}}+\pi_0 \label{p0_sw}\end{eqnarray}

Under these boundary conditions, the pressure is of order $\varepsilon^0$ instead of order $\varepsilon^{-2}$. In this case, from equations \eqref{ec_ui0_coefs} we obtain 
\begin{eqnarray}
&&\hspace*{-0.5cm}     
 \bar{u}_i^2 = -\varepsilon h \dfrac{A^1}{2 A^0}
\bar{u}_i^1+O(\varepsilon^2) \label{ui2_ui1}
 \quad i=1,2\end{eqnarray}
 
Boundary condition \eqref{Tn1_xi3_1} can be written (using
 \eqref{n1}) as follows:
\begin{eqnarray}
&&\left(\sigma_{ij}^{\varepsilon}
{v}^{\varepsilon}_{3j}\right)\cdot
{v}^{\varepsilon}_{3i}=-\pi^{\varepsilon}_1\|\vec{v}^{\varepsilon}_3\|^2
\textrm{ on } \xi_3=1 \label{Tn1_xi3_1s}
\end{eqnarray}

We use expressions  \eqref{sigmaij_des}, \eqref{v3} and \eqref{mod_v3}  to substitute $\sigma_{ij}$, vector $\vec{v}_3$ and its module into the above condition. Simplifying we yield
\begin{eqnarray}
&&{}\|\vec{a}_1 \times \vec{a}_2\|^2\left( -\sum_{k=0}^3 \bar{p}^k  +  \dfrac{2 \mu}{\varepsilon h} \sum_{k=1}^3  k   \bar{u}_3^k \right)-  \dfrac{2 \mu}{ h} \|\vec{a}_1 \times \vec{a}_2\|^2 \left( \dfrac{\partial h}{\partial
\xi_1} \sum_{k=1}^3  k   \bar{u}_1^k + \dfrac{\partial h}{\partial \xi_2} \sum_{k=1}^3  k   \bar{u}_2^k\right)\nonumber\\
&&{}+  2 \mu \|\vec{a}_1 \times \vec{a}_2\| \vec{a}_{3} \left[  \left(\vec{a}_1 \times \dfrac{\partial
\vec{a}_3}{\partial \xi_2}\right) +  \left(\dfrac{\partial
\vec{a}_3}{\partial \xi_1} \times
\vec{a}_2\right)\right]\sum_{k=1}^3  k   \bar{u}_3^k=-\pi_1\|\vec{a}_1 \times \vec{a}_2\|^2 + O(\varepsilon)
 \end{eqnarray}
 that can be rewritten, keeping in mind \eqref{p0_sw}, \eqref{p10}, \eqref{p20}, \eqref{p3e1}, \eqref{ui3_ui2}, \eqref{u_i^2_p0}, \eqref{ui1_WVui2}, \eqref{u_3^1}, \eqref{u32}, \eqref{u33}, in this way:
 \begin{eqnarray}
&&{}\pi_0 +  \dfrac{\mu}{\sqrt{A^0}} \textrm{div}\left(\sqrt{A^0} \left(\vec{W}- \vec{V}\right)\right)
 +  \dfrac{ \mu}{ h} \nabla h\cdot \left(\vec{W}- \vec{V}\right)=\pi_1 + O(\varepsilon)
 \end{eqnarray}

To impose boundary conditions  \eqref{Tai_xi3_0}-\eqref{Tvi_xi3_1} we introduce the friction force in the following typical way
\begin{equation}
    \vec{f}^{\varepsilon}_{R\alpha} = \rho_0 C_R^\varepsilon \|
\vec{{u}}^\varepsilon \| \vec{{u}}^\varepsilon \textrm{ on } \xi_3=\alpha, \quad (\alpha=0,1)\label{fReps}
\end{equation}
where $C_R^\varepsilon$ is a small constant. Let us assume that it is of order $\varepsilon$, that is, 
\begin{equation}
    C_R^{\varepsilon} =\varepsilon C^1_R\label{CR}
\end{equation}
so the friction forces can be written as follows
\begin{equation}
    \vec{f}^{\varepsilon}_{R\alpha} = \varepsilon  \vec{f}^{1}_{R\alpha} +O(\varepsilon^2) \textrm{ on } \xi_3=\alpha, \quad (\alpha=0,1)\label{fR}
\end{equation}

Now,  boundary conditions \eqref{Tai_xi3_0} on $\xi_3=0$ can be written  using
\eqref{sigmaij_des} and \eqref{n0}:
  \begin{eqnarray}
&&\hspace*{-0.5cm}\mu  \left[ \dfrac{\partial \bar{u}_3^0 }{\partial
\xi_1}+ \sum_{m=1}^2  
 \bar{u}_m^0 \left(\dfrac{\partial \vec{a}_{m}}{\partial
\xi_1} \cdot \vec{a}_{3}\right)+ \dfrac{1  }{\varepsilon h} \sum_{m=1}^2  \bar{u}_m^1 (\vec{a}_{m} \cdot  \vec{a}_{1}) \right]= {-} s_0  \varepsilon  \vec{f}^{1}_{R_0}  \cdot \vec{a}_{1}+O(\varepsilon^2) 
\label{desarrollo_cc_fR01}\\
&&\hspace*{-0.5cm} \mu \left[  \dfrac{\partial \bar{u}_3^0 }{\partial
\xi_2}+\sum_{m=1}^2  \bar{u}_m^0 \left(\dfrac{\partial \vec{a}_{m}}{\partial
\xi_2} \cdot \vec{a}_{3}\right)+ \dfrac{1 }{\varepsilon h}  \sum_{m=1}^3  \bar{u}_m^1 (\vec{ a}_{m} \cdot  \vec{a}_{2})\right]= {-} s_0  \varepsilon  \vec{f}^{1}_{R_0}  \cdot \vec{a}_{2} +O(\varepsilon^2)   \label{desarrollo_cc_fR02}
\end{eqnarray}

We multiply \eqref{desarrollo_cc_fR01} by $\alpha^0_i$ and \eqref{desarrollo_cc_fR02} by $\beta^0_i$, then we sum both equations and obtain
\begin{eqnarray}
&&\hspace*{-1.2cm}\bar{u}_i^1=-\varepsilon h \left[ \sum_{k=1}^2 \bar{u}_k^0 D^0_{ik}+ \alpha^0_i\dfrac{\partial \bar{u}_3^0}{\partial \xi_1}+ \beta^0_i\dfrac{\partial \bar{u}_3^0}{\partial \xi_2}+ \dfrac{s_0 \varepsilon}{\mu} (\alpha_i^0 \vec{f}^1_{R_0} \cdot \vec{a}_1+\beta_i^0 \vec{f}^1_{R_0} \cdot \vec{a}_2)\right] +O(\varepsilon^3)   \quad (i=1,2) \label{u_i^1sw}
\end{eqnarray} 

This implies that the terms $\bar{u}_i^2$ ($i=1,2,3$)  are of order $\varepsilon^2$ (\eqref{ui2_ui1},  \eqref{u32}) and, the terms $\bar{u}_i^3$ ($i=1,2,3$) are of order $\varepsilon^3$ (\eqref{ec_u1i}, \eqref{u33}). Considering this, \eqref{cc_u3_xi3_1b} yields
\begin{equation}
\bar{u}_3^1=\varepsilon \dfrac{\partial h}{\partial t}+O(\varepsilon^2) \label{u_3^1_dh_dt}
\end{equation}
and therefore \eqref{p0_sw} could be written as
\begin{eqnarray}
\bar{p}^0&=& \dfrac{2\mu }{ h} \dfrac{\partial h}{\partial t} +\pi_0+O(\varepsilon) \label{p0_swb}\end{eqnarray}

Boundary conditions \eqref{Tvi_xi3_1} can be rewritten taking into account in the first place \eqref{n1}, \eqref{v1},  \eqref{v2}, \eqref{sigmaij_des}, \eqref{v3}, \eqref{mod_v3} and \eqref{fR}. Then, we have grouped and divided by $ \|\vec{a}_1 \times \vec{a}_2\| $. And, considering that $\bar{u}_i^3$, $i=1,2$ are of order $\varepsilon^3$ and $\bar{u}_i^2$ are of order $\varepsilon^2$ we get
\begin{eqnarray}
&&\hspace*{-0.5cm} \dfrac{1}{\varepsilon h}  \sum_{k=1}^2 k\left[    \bar{u}_1^k  (\vec{a}_{1} \cdot \vec{a}_i )  + \bar{u}_2^k  (\vec{a}_{2} \cdot \vec{a}_i )\right]+\dfrac{1}{ h}  
\dfrac{\partial h}{\partial \xi_i}     \bar{u}_3^1  + \sum_{m=1}^2 
 \bar{u}_m^0   \left( \dfrac{\partial \vec{a}_{m}}{\partial
\xi_i} \cdot\vec{a}_3\right)+ \dfrac{\partial \bar{u}_3^0 }{\partial
\xi_i} 
 \nonumber\\
&&{}+     \bar{u}_1^1 \left[ \left(   \vec{a}_{1} \cdot
\dfrac{\partial \vec{a}_3}{\partial \xi_i} \right)+ \dfrac{ (\vec{a}_{1} \cdot \vec{a}_i ) }{\|\vec{a}_1 \times \vec{a}_2\| } I \right] + \bar{u}_2^1 \left[  \left( \vec{a}_{2} \cdot
\dfrac{\partial \vec{a}_3}{\partial \xi_i} \right) + \dfrac{(\vec{a}_{2} \cdot \vec{a}_i )}{\|\vec{a}_1 \times \vec{a}_2\| } I\right] \nonumber\\
&&{}+   \dfrac{\partial \bar{u}_3^1 }{\partial
\xi_i} + \sum_{m=1}^2 \bar{u}_m^1 \left(\dfrac{\partial \vec{a}_{m}}{\partial
\xi_i} \cdot \vec{a}_3 \right)
 - \varepsilon \sum_{m=1}^2 \sum_{l=1}^2 
\dfrac{\partial( \bar{u}_m^0  \vec{a}_{m})}{\partial \xi_l} \cdot \vec{a}_i \left( 
\left(\alpha_l^0\dfrac{\partial h}{\partial
\xi_1}  + \beta_l^0 \dfrac{\partial h}{\partial \xi_2} \right)\right)+\varepsilon \bar{u}_3^0 H^0_{3i3}\nonumber\\
&&{}- \varepsilon \left( \dfrac{\partial \bar{u}_1^0 }{\partial
\xi_i}  
 \dfrac{\partial h}{\partial
\xi_1}+ \dfrac{\partial \bar{u}_2^0 }{\partial
\xi_i}  
 \dfrac{\partial h}{\partial
\xi_2}\right) + \dfrac{\varepsilon}{\|\vec{a}_1 \times \vec{a}_2\| }\sum_{m=1}^3 \bar{u}_m^0  \dfrac{\partial \vec{a}_{m}}{\partial
\xi_i}\cdot 
\vec{\eta}(h)
\nonumber\\
&&{} +\dfrac{\varepsilon h }{\|\vec{a}_1 \times \vec{a}_2\| } \dfrac{\partial \bar{u}_3^0 }{\partial
\xi_i} I
=\dfrac{s_0}{\mu} \varepsilon \vec{f}^{1}_{R_1}\cdot
 \vec{a}_i+O(\varepsilon^2), \quad (i=1,2)
\end{eqnarray}
where $H^0_{3i3}$, $I$ and $\vec{\eta}(h)$ are given by \eqref{H}, \eqref{I} and \eqref{eta}.

Next, we substitute $\bar{u}_i^k$, $i=1,2,3, k=0,1$ using expressions \eqref{ui0_Vi}, \eqref{u30}, \eqref{u_i^1sw} and \eqref{u_3^1_dh_dt}
\begin{eqnarray}
&&\hspace*{-0.5cm} 2\left[    \bar{u}_1^2   (\vec{a}_{1} \cdot \vec{a}_i )   + \bar{u}_2^2 (\vec{a}_{2} \cdot \vec{a}_i )\right]= \varepsilon^2 h \sum_{k=1}^2 \sum_{l=1}^2 
\dfrac{\partial( V_k  \vec{a}_{k})}{\partial \xi_l} \cdot \vec{a}_i 
\left(\alpha_l^0\dfrac{\partial h}{\partial
\xi_1}  + \beta_l^0 \dfrac{\partial h}{\partial \xi_2} \right) \nonumber\\
&&{}+ \varepsilon^2 h \left( \dfrac{\partial V_1 }{\partial
\xi_i}  
 \dfrac{\partial h}{\partial
\xi_1}+ \dfrac{\partial V_2 }{\partial
\xi_i}  
 \dfrac{\partial h}{\partial
\xi_2}\right) + \dfrac{\varepsilon^2 h }{\|\vec{a}_1 \times \vec{a}_2\| } \sum_{k=1}^2 V_k \left[hI\left( \vec{a}_3 \cdot \dfrac{\partial \vec{a}_k}{\partial \xi_i}\right) -\dfrac{\partial \vec{a}_{k}}{\partial
\xi_i}\cdot 
\vec{\eta}(h)\right] \nonumber\\
&&{}
-\varepsilon^2 
\dfrac{\partial h}{\partial \xi_i}     \dfrac{\partial h}{\partial t}  - \varepsilon^2 h  \dfrac{\partial^2 h}{\partial
\xi_i \partial t}  
-\varepsilon^2 h \left( \dfrac{\partial X}{\partial t} \cdot \vec{a}_3\right) \left[H^0_{3i3} + \dfrac{1}{\|\vec{a}_1 \times \vec{a}_2\| }\left(  \dfrac{\partial \vec{a}_{3}}{\partial
\xi_i}\cdot 
\vec{\eta}(h)\right)\right]\nonumber\\
&&{}
+\varepsilon^2 h\dfrac{s_0}{\mu}\left( \vec{f}^{1}_{R_1}
 +\vec{f}^{1}_{R_0}\right) \cdot
 \vec{a}_i +O(\varepsilon^3), \quad (i=1,2) \label{ui21_desde_cc}
\end{eqnarray}

We sum equation \eqref{ui21_desde_cc} ($i=1$) multiplied by $\alpha_j^0$ and equation \eqref{ui21_desde_cc} ($i=2$)  multiplied by $\beta_j^0$, $j=1,2$. In this way, we obtain the following expressions for $\bar{u}_i^2$, $i=1,2$
\begin{eqnarray}
&&\hspace*{-0.5cm}     \bar{u}_i^2 = \dfrac{\varepsilon^2 h}{2} \left\{ \sum_{l=1}^2 \left[-J^{0,0}_{3l}
\dfrac{\partial V_i }{\partial \xi_l}
+ J^{0,0}_{li} \sum_{k=1}^2 \dfrac{\partial V_k }{\partial
\xi_l}  
 \dfrac{\partial h}{\partial
\xi_k} -J^{0,0}_{3l} \sum_{k=1}^2 V_k H^0_{ilk} \right] \right.\nonumber\\
&&{} + \dfrac{1 }{\sqrt{A^0} } \sum_{k=1}^2 V_k \left[hI D^0_{ik} - \alpha^0_i \left(\dfrac{\partial \vec{a}_{k}}{\partial
\xi_1}\cdot 
\vec{\eta}(h)\right) - \beta^0_i \left(\dfrac{\partial \vec{a}_{k}}{\partial
\xi_2}\cdot 
\vec{\eta}(h)\right)\right]
+\dfrac{ J^{0,0}_{3i}}{h}  \dfrac{\partial h}{\partial t} \nonumber\\
&&\left.{} -  \sum_{k=1}^2 J^{0,0}_{ki} \dfrac{\partial^2 h}{\partial
\xi_k \partial t}  
-\left( \dfrac{\partial X}{\partial t} \cdot \vec{a}_3\right) \left[ \sum_{l=1}^2 J^{0,0}_{3l}
H^0_{il3} 
 + \dfrac{J^{0,0}_{li}}{\sqrt{A^0} } \left(  \dfrac{\partial \vec{a}_{3}}{\partial
\xi_l}\cdot 
\vec{\eta}(h)\right)\right]\right\}\nonumber\\
&&{}
+\dfrac{\varepsilon^2 h}{2}\dfrac{s_0}{\mu}\left[ \alpha^0_i\left( \vec{f}^{1}_{R_1}
 +\vec{f}^{1}_{R_0}\right) \cdot
 \vec{a}_1 + \beta^0_i\left( \vec{f}^{1}_{R_1}
 +\vec{f}^{1}_{R_0}\right) \cdot
 \vec{a}_2 \right] +O(\varepsilon^3), \quad (i=1,2)\label{u12_desde_cc}
\end{eqnarray}

Now, we can rewrite equations \eqref{ec_ui0_coefs} taking into account  \eqref{ui0_Vi} and \eqref{u30} and, substituting  $\bar{p}^0$, $\bar{u}_i^1$  and $\bar{u}_i^2$ using  \eqref{p0_swb}, \eqref{u_i^1sw} and \eqref{u12_desde_cc} 
 \begin{eqnarray}
&&\hspace*{-0.5cm} \dfrac{ \partial V_i}{\partial t}+\sum_{l=1}^2  \dfrac{ \partial V_i}{\partial \xi_l} (V_l -C^0_l)+ \sum_{k=1}^2
 V_k  \left(  Q^0_{ik} + \sum_{l=1}^2 V_l H^0_{ilk}+ H^0_{ik3} \left( \dfrac{\partial \vec{X}}{\partial t} \cdot \vec{a}_3\right)  \right)
 \nonumber\\
&&{}
 =-\dfrac{1}{\rho_0} \sum_{l=1}^2 \dfrac{
\partial  }{\partial \xi_l}\left(  \dfrac{2\mu }{ h}  \dfrac{\partial h}{\partial t}+\pi_0 \right) J^{0,0}_{il}  + \nu \left\{ \sum_{m=1}^2 \sum_{l=1}^2 \dfrac{
\partial^2 V_i}{\partial \xi_l \partial \xi_m} J^{0,0}_{lm} + \sum_{k=1}^2 \sum_{l=1}^2 \dfrac{
\partial V_k}{ \partial \xi_l} L^0_{kli}\right. \nonumber\\
&&\left.{} +  \sum_{k=1}^2  V_k {S}^0_{ik}+  \sum_{l=1}^2 \dfrac{
\partial }{ \partial \xi_l}\left( \dfrac{\partial \vec{X}}{\partial t} \cdot \vec{a}_3\right) L^0_{3li}
+  \left( \dfrac{\partial \vec{X}}{\partial t} \cdot \vec{a}_3\right) {S}^0_{i3}\right. \nonumber\\
&&\left.{}-  \dfrac{A^1}{A^0} \sum_{k=1}^2
\left[ V^0_k D^0_{ik}+ J^{0,0}_{ik}\dfrac{\partial }{\partial \xi_k}\left( \dfrac{\partial \vec{X}}{\partial t} \cdot \vec{a}_3\right)\right]  + 
\dfrac{1}{h}\sum_{l=1}^2 \left[ -J^{0,0}_{3l}
\dfrac{\partial V_i }{\partial \xi_l}
+  J^{0,0}_{li} \sum_{k=1}^2 \dfrac{\partial V_k }{\partial
\xi_l}  
 \dfrac{\partial h}{\partial
\xi_k}\right]\right.\nonumber\\
&&{}- \dfrac{1}{h} \sum_{k=1}^2 V_k\sum_{l=1}^2 
H^0_{ilk} 
J^{0,0}_{3l} + \dfrac{1}{h \sqrt{ A^0} } \sum_{k=1}^2 V_k \left[hI D^0_{ik}- \alpha^0_i \left(\dfrac{\partial \vec{a}_{k}}{\partial
\xi_1}\cdot 
\vec{\eta}(h)\right) - \beta^0_i \left(\dfrac{\partial \vec{a}_{k}}{\partial
\xi_2}\cdot 
\vec{\eta}(h)\right)\right]\nonumber\\
&&{}
+\dfrac{J^{0,0}_{3i}}{h^2}
   \dfrac{\partial h}{\partial t}  - \dfrac{1}{h}\sum_{k=1}^2 J^{0,0}_{ki} \dfrac{\partial^2 h}{\partial
\xi_k \partial t}  
-\dfrac{1}{h} \left( \dfrac{\partial X}{\partial t} \cdot \vec{a}_3\right) \left[ \sum_{l=1}^2 J^{0,0}_{3l}
H^0_{il3} 
 + \dfrac{J^{0,0}_{li}}{\sqrt{A^0} } \left(  \dfrac{\partial \vec{a}_{3}}{\partial
\xi_l}\cdot 
\vec{\eta}(h)\right)  \right]\nonumber\\
&&\left.{}
+\dfrac{1}{ h}\dfrac{s_0}{\mu}\left[ \alpha^0_i\left( \vec{f}^{1}_{R_1}
 +\vec{f}^{1}_{R_0}\right) \cdot
 \vec{a}_1 + \beta^0_i\left( \vec{f}^{1}_{R_1}
 +\vec{f}^{1}_{R_0}\right) \cdot
 \vec{a}_2 \right]  
 \right\}   + {\bar{f}_i^0} - 
\left( \dfrac{\partial \vec{X}}{\partial t} \cdot \vec{a}_3\right)  Q^0_{i3}+O(\varepsilon)\nonumber\\
&&(i=1,2)\end{eqnarray}

That can be written in a more compact form
 \begin{eqnarray}
&&\hspace*{-0.5cm} \dfrac{ \partial V_i}{\partial t}+\sum_{l=1}^2  \dfrac{ \partial V_i}{\partial \xi_l} (V_l -C^0_l)+ \sum_{k=1}^2
 V_k  \left(  R^0_{ik} + \sum_{l=1}^2 V_l H^0_{ilk} \right) =-\dfrac{1}{\rho_0} \sum_{l=1}^2 \dfrac{
\partial  \pi_0  }{\partial \xi_l}J^{0,0}_{il}
 \nonumber\\
&&{}
  + \nu \left\{ \sum_{m=1}^2 \sum_{l=1}^2 \dfrac{
\partial^2 V_i}{\partial \xi_l \partial \xi_m} J^{0,0}_{lm} + \sum_{k=1}^2 \sum_{l=1}^2 \dfrac{
\partial V_k}{ \partial \xi_l} ( L^0_{kli} +\psi(h)^0_{ikl} ) +  \sum_{k=1}^2  V_k \bar{S}^0_{ik} + \kappa(h)^0_i \right\}  \nonumber\\
&&{} - 
\left( \dfrac{\partial \vec{X}}{\partial t} \cdot \vec{a}_3\right)  Q^0_{i3} +\bar{F}^0_i(h)+O(\varepsilon), \quad (i=1,2) \label{ec_Vi0-new} \end{eqnarray}
where coefficients $C^0_l$, $R^0_{ik}$, $H^0_{ilk}$, $J^{0,0}_{lm}$, $L^0_{kli}$, $\psi(h)^0_{ikl}$, $\bar{S}^0_{ik}$, ${\kappa}(h)^0_i$, $\bar{F}^0_i(h)$ and $Q^0_{i3}$ are defined in appendix \ref{ApendiceA} (in \eqref{C}, \eqref{R}, \eqref{H}, \eqref{L}, \eqref{psi}, \eqref{barS}, \eqref{kappa}, \eqref{Fbar} and \eqref{Q3} respectively).

Keeping in mind that $\pi_0 = \pi_0^0 + O(\varepsilon)$, \eqref{F}-\eqref{Fbar} and \eqref{barSPchi}, we deduce that equations \eqref{ec_Vi0-new} are equivalent to equations \eqref{ec_Vi0-v2}.

\section{New model} \label{NewModel}

As we have just seen in sections \ref{NewAssumptions} and \ref{BoundaryConditions}, hypotheses \eqref{ui_pol_xi3}-\eqref{f_pol_xi3} allow us to obtain a two-dimensional model formed by equations \eqref{u_3^1}, \eqref{ec_ui0_coefs}-\eqref{ui_3}, which we shall call new model from now on, whose asymptotic behavior, when $\varepsilon$ tends to zero, is the same as the Navier-Stokes equations. 

In fact, we have justified that, under the assumptions about the boundary conditions made in subsection \ref{subseccion-4-1}, the solution of the new model approaches the solution of model \eqref{lubric} as $\varepsilon$ tends to zero, just as in the previous work \cite{RodTabJMAA2021}, where we showed that the solution of the Navier-Stokes equations
approaches the solution of \eqref{Reynolds_gen}.  And, we have also seen that, under the assumptions about the boundary conditions shown in subsection \ref{subseccion-4-2}, the solution of the the new model tends to the solution of \eqref{ec_Vi0-new}, as it happened in our prior article \cite{RodTabJMAA2021} with the solution of the Navier-Stokes equations (see \eqref{ec_Vi0-v2}).

Examining the new model, we observe that the equations can be divided into  two groups: a first group, including equations \eqref{ec_ui0_coefs}, \eqref{ec_u1i}, \eqref{ui_2}, \eqref{ui_3} and \eqref{div_i_dr_alfa_beta_pol_xi3_3s}, that must be solved to obtain the terms  $\bar{u}_1^0$, $\bar{u}_2^0$, $\bar{p}^0$, $\bar{u}_1^1$, $\bar{u}_2^1$, $\bar{u}_1^2$, $\bar{u}_2^2$, $\bar{u}_1^3$ and $\bar{u}_2^3$, and a second group, including equations \eqref{u30}, \eqref{u_3^1}, \eqref{u32}, \eqref{u33}, \eqref{ec_ns_u30_Res}, \eqref{ec_u31} and \eqref{u3_2}, that allow us to eliminate the terms $\bar{u}_3^0$, $\bar{u}_3^1$, $\bar{u}_3^2$, $\bar{u}_3^3$, $\bar{p}^1$, $\bar{p}^2$, $\bar{p}^3$ from equations \eqref{ec_ui0_coefs}, \eqref{ec_u1i}, \eqref{ui_2}, \eqref{ui_3} and  \eqref{div_i_dr_alfa_beta_pol_xi3_3s}. Once the aforementioned elimination has been carried out, we can solve the first group of equations to compute $\bar{u}_i^k$ $(k=0,1,2,3, i=1,2)$ and $\bar{p}^0$, and use the second group of equations to obtain $\bar{u}_3^k$ ($k=0,1,2,3$) and $\bar{p}^j$ ($j=1,2,3$). Boundary and initial conditions must be added to this system of equations too.

\section{Conclusions}

In this article, we propose a two-dimensional flow model of a viscous fluid between two very close moving surfaces and, we show that its asymptotic behavior, when the distance between the two surfaces tends to zero, is the same as that previously obtained \cite{RodTabJMAA2021} for the Navier-Stokes equations. Depending on the boundary conditions of the problem, the solutions of both models converge to the solutions of two different limit problems: in the first case, when slip velocity boundary conditions on the upper and lower surfaces are imposed, a lubrication model is derived (see subsection \ref{subseccion-4-1}); if, instead, the tractions and friction forces are known on both bound surfaces, a shallow water model is obtained (see subsection \ref{subseccion-4-2}).

As it is well known, numerical solution of three-dimensional Navier-Stokes equations requires large computational resources, and solving these equations in such a thin domain presents even more numerical problems, while solving the new two-dimensional model presented here is much easier. 

 On the other hand, as we have already mentioned, the new model has the same asymptotic behavior as the Navier-Stokes equations, so, in a certain sense, it encompasses the two limit models presented in subsections \ref{subseccion-4-1} and \ref{subseccion-4-2}.
 
For all the above reasons, the new model proposed in this article can be considered a good option for calculating viscous fluid flow between two very close moving surfaces, without the need to decide a priori whether the flow is typical of a lubrication problem or it is of shallow water type, and without the enormous computational effort that would be required to solve the Navier-Stokes equations in such a thin domain.
 
It remains for us to carry out numerical simulations to compare the accuracy and the computation time needed to solve the new and the classic models that we have mentioned. We hope to present results in this regard very soon.

\section*{Acknowledgements}
This work has been partially supported by  the European Union's Horizon 2020 Research and Innovation Programme,
under the Marie Sklodowska-Curie Grant Agreement No 823731 CONMECH. 

\appendix

\section{Coefficients definition} \label{ApendiceA}

In this appendix, we introduce some coefficients that depend either only on the lower bound surface parametrization, $\vec{X}$ or on both the parametrization and the gap $h$. We will use these coefficients throughout this article.

In the first place, the coefficients of the first and second fundamental forms of the surface parametrized by $\vec{X}$ have been denoted by $E, F, G$ and $e,f,g$, respectively:
\begin{eqnarray} E&=&\vec{a}_1 \cdot\vec{a}_1 \label{coef-E} \\
 F&=&\vec{a}_1 \cdot\vec{a}_2 \label{coef-F} \\
 G&=&\vec{a}_2 \cdot\vec{a}_2 \label{coef-G}\\
e&=& -\vec{a}_1 \cdot \dfrac{\partial \vec{a}_{3}}{\partial
\xi_1} = \vec{a}_3 \cdot \dfrac{\partial \vec{a}_{1}}{\partial
\xi_1} \label{coef-e} \\
f&=&-
\vec{a}_1 \cdot \dfrac{\partial \vec{a}_{3}}{\partial
\xi_2}=-\vec{a}_2 \cdot\dfrac{\partial \vec{a}_{3}}{\partial
\xi_1}= \vec{a}_3 \cdot \dfrac{\partial \vec{a}_{1}}{\partial
\xi_2}= \vec{a}_3 \cdot \dfrac{\partial \vec{a}_{2}}{\partial
\xi_1} \label{coef-f} \\
g&=&-\vec{a}_2 \cdot\dfrac{\partial \vec{a}_{3}}{\partial
\xi_2}= \vec{a}_3 \cdot \dfrac{\partial \vec{a}_{2}}{\partial
\xi_2} \label{coef-g} 
\end{eqnarray}
and, from them, we define:

\begin{eqnarray}
A^0&=&\|\vec{a}_1\|^2 \|\vec{a}_2\|^2-\left( \vec{a}_1
\cdot \vec{a}_{2}\right)^2 =EG-F^2=\|\vec{a}_1 \times \vec{a}_2\|^2 \label{A0}  \\
A^1&=&\|\vec{a}_2\|^2 \left(\vec{a}_1 \cdot\dfrac{\partial
\vec{a}_{3}}{\partial \xi_1}\right) + \|\vec{a}_1\|^2
\left(\vec{a}_2 \cdot \dfrac{\partial \vec{a}_{3}}{\partial
\xi_2}\right) \nonumber \\
&&{} - \left( \vec{a}_1 \cdot
\vec{a}_{2}\right)\left(\vec{a}_1 \cdot\dfrac{\partial
\vec{a}_{3}}{\partial \xi_2}+\vec{a}_2 \cdot \dfrac{\partial
\vec{a}_{3}}{\partial \xi_1}\right) = -eG-gE+2fF \label{A^1} \\ 
A^2&=&\left(\vec{a}_1 \cdot\dfrac{\partial
\vec{a}_{3}}{\partial \xi_1} \right)\left(\vec{a}_2 \cdot
\dfrac{\partial \vec{a}_{3}}{\partial \xi_2} \right) -\left(
\vec{a}_1 \cdot\dfrac{\partial \vec{a}_{3}}{\partial \xi_2}\right)
\left( \vec{a}_2 \cdot \dfrac{\partial \vec{a}_{3}}{\partial \xi_1}
\right)=eg-f^2\label{A2}\\
M&=&\begin{pmatrix}G &
-F\\
-F &
E\end{pmatrix}\label{M}
\end{eqnarray}

The following coefficients are involved in the change of variable defined in section \ref{PreviousModels}:
\begin{eqnarray}
\alpha_i
&=&
\alpha_i^0+\varepsilon \xi_3 h \alpha_i^1 +\varepsilon^2 \xi_3^2h^2
\alpha_i^2+\cdots,
\quad (i=1,2) \label{alfaides}\\
\alpha_3 &=&\dfrac{ \xi_3 }{ h}(\alpha_3^0+\varepsilon \xi_3 h \alpha_3^1
+\varepsilon^2 \xi_3^2h^2 \alpha_3^2+\cdots), \label{alfa3des}\\
 \beta_i &=&\beta_i^0+\varepsilon \xi_3 h \beta_i^1 +\varepsilon^2 \xi_3^2h^2
\beta_i^2+\cdots, \quad (i=1,2) \label{betaides}\\
\beta_3&=&\dfrac{ \xi_3 }{ h}(\beta_3^0+\varepsilon \xi_3 h \beta_3^1
+\varepsilon^2 \xi_3^2h^2 \beta_3^2+\cdots), \label{beta3des} \\
\gamma_3&=&\dfrac{1}{\varepsilon h}, \quad \gamma_1 = \gamma_2 = 0, \label{gammades} 
\end{eqnarray}
where 
\begin{eqnarray}
&&\alpha_1^0=\dfrac{\|\vec{a}_2\|^2}{A^0}=\dfrac{G}{EG-F^2}\label{alfa10}\\
&&\alpha_1^1= \dfrac{\vec{a}_2\cdot \dfrac{\partial
\vec{a}_{3}}{\partial \xi_2}- \alpha_1^0 A^1}{ A^0}=- \dfrac{g +
\alpha_1^0 A^1}{ A^0}\label{alfa11}
\\
&&\alpha_1^n= -\dfrac{\alpha_1^{n-2} A^2+\alpha_1^{n-1} A^1}{  A^0},
\quad n\geq 2
\label{alfa1n}\\
&&\alpha_2^0=\beta_1^0=-\dfrac{ \vec{a}_2\cdot \vec{a}_{1}}{ A^0} =-\dfrac{ F}{ A^0}  \label{alfa20_beta10}\\
&&\alpha_2^1=\beta_1^1 =- \dfrac{\vec{a}_2\cdot \dfrac{\partial
\vec{a}_{3}}{\partial \xi_1} +\alpha_2^0 A^1}{ A^0}=
\dfrac{f- \alpha_2^0 A^1}{ A^0} \label{alfa21_beta11}\\
&&\alpha_2^n = \beta_1^n= -\dfrac{ \alpha_2^{n-2} A^2 +
\alpha_2^{n-1} A^1}{ A^0},
\quad n\geq 2  \label{alfa2n_beta1n}\\
&&\alpha_3^0=\dfrac{\dfrac{\partial h}{\partial \xi_2} \vec{a}_{1}
\cdot \vec{a}_2 - \dfrac{\partial h}{\partial \xi_1} \|\vec{a}_2\|^2
}{A^0} =-\alpha_1^0 \dfrac{\partial h}{\partial \xi_1}  -
\alpha_2^0\dfrac{\partial h}{\partial \xi_2}    \label{alfa30}\\
&&\alpha_3^1= \dfrac{\vec{a}_2\cdot\left[\dfrac{\partial h}{\partial
\xi_2}
  \dfrac{\partial
\vec{a}_{3}}{\partial \xi_1} - \dfrac{\partial h}{\partial \xi_1}
\dfrac{\partial \vec{a}_{3}}{\partial \xi_2}\right]-\alpha_3 ^0 A^1
} { A^0} =-\alpha_1^1 \dfrac{\partial h}{\partial \xi_1}  -
\alpha_2^1\dfrac{\partial h}{\partial \xi_2}
\label{alfa31}\\
 && \alpha_3^n=-\dfrac{\alpha_3^{n-2} A^2 + \alpha_3^{n-1} A^1} {A^0},  \quad n\geq 2 \label{alfa3n}
\\
&&\beta_2^0=\dfrac{ \|\vec{a}_{1}\|^2}{ A^0} =\dfrac{ E}{ A^0}  \label{beta20}\\
&&\beta_2^1 = \dfrac{\vec{a}_1\cdot \dfrac{\partial
\vec{a}_{3}}{\partial \xi_1} -\beta_2^0 A^1}{ A^0} =-
\dfrac{ e +\beta_2^0 A^1}{ A^0}\label{beta21}\\
&&\beta_2^n =-\dfrac{ \beta_2^{n-2} A^2 + \beta_2^{n-1} A^1}{ A^0},  \quad n\geq 2   \label{beta22}\\
&&\beta_3^0=\dfrac{\dfrac{\partial h}{\partial \xi_1} \vec{a}_{1}
\cdot \vec{a}_2 - \dfrac{\partial h}{\partial \xi_2}\|\vec{a}_1\|^2
}{A^0} =-\beta_1^0 \dfrac{\partial h}{\partial \xi_1}  -
\beta_2^0\dfrac{\partial h}{\partial \xi_2} \label{beta30}\\
&&\beta_3^1= \dfrac{\dfrac{\partial h}{\partial \xi_1} \left(
\vec{a}_1 \cdot\dfrac{\partial \vec{a}_{3}}{\partial \xi_2}\right) -
\dfrac{\partial h}{\partial \xi_2} \left(  \vec{a}_1
\cdot\dfrac{\partial \vec{a}_{3}}{\partial \xi_1}\right)-\beta_3 ^0
A^1 } { A^0} = -\beta_1^1 \dfrac{\partial h}{\partial \xi_1}  -
\beta_2^1\dfrac{\partial h}{\partial \xi_2}
\label{beta31}\\
 && \beta_3^n=-\dfrac{\beta_3^{n-2} A^2 + \beta_3^{n-1} A^1} {A^0},  \quad n\geq 2
 \label{beta3n}
\end{eqnarray}

The next set of coefficients
depend only on the parametrization $\vec{X}$, except for those that include $\alpha_3^j$ and $\beta_3^j$ that also depend on $h$: 
\begin{eqnarray}
B^j_{lk}&=& \alpha^j_l (\vec{a}_1 \cdot \vec{a}_k) +
 \beta^j_l (\vec{a}_2 \cdot \vec{a}_k), \quad (j=0,1,2;\quad l=1,2,3;\quad k=1,2) \label{B}\\
C^0_l&=& \alpha_l^0 \left( \vec{a}_1 \cdot \dfrac{\partial \vec{X}}{\partial t}\right) + \beta_l^0 \left( \vec{a}_2 \cdot \dfrac{\partial \vec{X}}{\partial t}\right) \quad (l=1,2,3) \label{C}\\
C^{i,j}_l&=& \alpha_l^i \left( \vec{a}_1 \cdot \dfrac{\partial \vec{X}}{\partial t}\right) + \beta_l^i \left( \vec{a}_2 \cdot \dfrac{\partial \vec{X}}{\partial t}\right)\nonumber\\
&+&\alpha_l^j \left( \vec{a}_1 \cdot \dfrac{\partial \vec{a}_3}{\partial t}\right) +
\beta_l^j \left( \vec{a}_2 \cdot \dfrac{\partial \vec{a}_3}{\partial t}\right) \quad (l=1,2,3; \quad i=1,2; \quad j=0,1,2) \label{Cij}\\
{D^j_{ik}}&=& \alpha_i^j\left(\vec{a}_3\cdot \dfrac{\partial \vec{a}_k}{\partial \xi_1}\right)+\beta_i^j\left(\vec{a}_3 \cdot \dfrac{\partial \vec{a}_k}{\partial \xi_2}\right)\quad (i=1,2,3; \quad k=1,2; \quad j=0,1)\label{D}\\
H^j_{ilk}&=& \alpha_i^j\left(\vec{a}_1\cdot \dfrac{\partial \vec{a}_k}{\partial \xi_l}\right)+\beta_i^j\left(\vec{a}_2 \cdot \dfrac{\partial \vec{a}_k}{\partial \xi_l}\right)\quad (l=1,2; \quad i,k=1,2,3; \quad j=0,1,2)\label{H}\\
I&=&\left(\vec{a}_1 \times \dfrac{\partial
\vec{a}_3}{\partial \xi_2}\right)\cdot \vec{a}_3 + \left(\dfrac{\partial
\vec{a}_3}{\partial \xi_1} \times
\vec{a}_2\right)\cdot \vec{a}_3 \label{I}\\
 J^{i,j}_{lm}&=& \alpha^i_l B^j_{m1} +\beta^i_l B^j_{m2} \quad (l, m=1,2,3; \quad i,j=0,1,2) \label{J}\\
 {{K}^{j,i}_{l}} &=&\sum_{m=1}^2
 \dfrac{\partial \alpha^j_l}{\partial \xi_m}B^i_{m1}+ \dfrac{\partial \beta^j_l}{\partial \xi_m}  B^i_{m2}+ \alpha^j_l H^i_{mm1} + \beta^j_l H^i_{mm2}, \quad (l=1,2,3; \quad i,j=0,1)
 \label{Kji}\\
{L}^0_{kli} &=&  {{K}^{0,0}_{l}} \delta_{ki}+ 2 \sum_{m=1}^2 H^0_{imk} J^{0,0}_{lm} \quad (i,l=1,2; \quad k=1,2,3) \label{L}\\
{L}^0_{kl3} &=& {{K}^{0,0}_{l}} \delta_{k3} +
2 D^0_{lk},  \quad (l=1,2; \quad k=1,2,3)  \label{L3}\\
P^0_{ik} &=&  \dfrac{I \sqrt{A^0} -A^1}{A^0} \left[ \alpha_i^0\left( \dfrac{\partial \vec{a}_k}{\partial \xi_1} \cdot \vec{a}_3\right) +\beta_i^0 \left( \dfrac{\partial \vec{a}_k}{\partial \xi_2} \cdot \vec{a}_3\right) \right]\nonumber\\
&+& \sum_{l=1}^2 \left[\sum_{m=1}^2  \left( \alpha^0_i \left( \vec{a}_1 \cdot \dfrac{\partial^2 \vec{a}_k}{\partial \xi_l \partial \xi_m}\right) +  \beta^0_i \left( \vec{a}_2 \cdot \dfrac{\partial^2 \vec{a}_k}{\partial \xi_l \partial \xi_m}\right)\right)J^{0,0}_{lm}+K^{0,0}_{l} H^0_{ilk}\right]\nonumber\\
&-&\dfrac{1}{\sqrt{A^0}}\left[ \left(\vec{a}_1 \times \dfrac{\partial
\vec{a}_3}{\partial \xi_2}\right)+ \left(\dfrac{\partial
\vec{a}_3}{\partial \xi_1} \times
\vec{a}_2\right)\right] \cdot \left( \alpha^0_i \dfrac{\partial \vec{a}_k}{\partial \xi_1} + \beta^0_i\dfrac{\partial \vec{a}_k}{\partial \xi_2} \right), \quad (i,k=1,2)\label{P} \\
Q^0_{ik}&=&\alpha_i^0 \left( \vec{a}_1 \cdot \dfrac{\partial \vec{a}_k}{\partial t}\right) + \beta_i^0 \left( \vec{a}_2 \cdot \dfrac{\partial \vec{a}_k}{\partial t}\right)-\sum_{l=1}^2 H^0_{ilk}C^0_l, \quad (i=1,2; \quad k=1,2,3) \label{Q} \\
{Q^0_{3k}}&=&
\left(\vec{a}_3 \cdot \dfrac{
\partial \vec{a}_{k}}{\partial t}\right) -\sum_{l=1}^2
 \left(\vec{a}_3 \cdot \dfrac{
\partial \vec{a}_{k}}{\partial \xi_l}\right) C^0_l \quad (k=1,2)\label{Q3}\\
R^0_{ik} &=&  Q^0_{ik}+ H^0_{ik3} \left( \dfrac{\partial \vec{X}}{\partial t} \cdot \vec{a}_3\right) \quad (i=1,2; \quad k=1,2) \label{R} 
\\
{{S}^0_{ik}} &=& \sum_{l=1}^2 \left[  \sum_{m=1}^2 \left( \alpha^0_i \left( \vec{a}_1 \cdot \dfrac{\partial^2 \vec{a}_k}{\partial \xi_l \partial \xi_m}\right) +  \beta^0_i \left( \vec{a}_2 \cdot \dfrac{\partial^2 \vec{a}_k}{\partial \xi_l \partial \xi_m}\right)\right)J^{0,0}_{lm} +  K^{0,0}_{l} H^0_{ilk}  \right]\nonumber\\
&&\quad
(i=1,2; \quad k=1,2,3)\label{S} \\
{{S}^0_{3k}} &=&
 \sum_{l=1}^2 \left[   \sum_{m=1}^2\left( \vec{a}_3\cdot\dfrac{
\partial^2 \vec{a}_{k} }{\partial \xi_l \partial \xi_m}\right)
J^{0,0}_{lm}+  \left( \vec{a}_3\cdot\dfrac{
\partial \vec{a}_{k}  }{\partial \xi_l}\right) K^{0,0}_{l}\right]\quad
( k=1,2,3)\label{S3}
\end{eqnarray} 

And, finally, we have the coefficientes that depend explicitly on function $h$ in addition to depending on the parametrization $\vec{X}$:
\begin{eqnarray}
F^0_i(h)&=&\int_0^1 f^0_i\, d\xi_3 +\dfrac{s_0}{\rho_0 h}(\vec{f}^1_{R_1} + \vec{f}^1_{R_0}) \cdot \left( \alpha^0_i \vec{a}_1 + \beta^0_i \vec{a}_2 \right) \quad (i=1,2) \label{F} \\
\bar{F}^0_i(h)&=&\bar{f}^0_i +\dfrac{s_0}{\rho_0 h}(\vec{f}^1_{R_1} + \vec{f}^1_{R_0}) \cdot \left( \alpha^0_i \vec{a}_1 + \beta^0_i \vec{a}_2 \right)=F^0_i(h)+O(\varepsilon) \quad (i=1,2) \label{Fbar} \\
{L}^{1,0}_{kli}(h) &=& {L}^0_{kli} + \dfrac{2}{h} J^{0,0}_{3l}\delta_{ki} \quad (l=1,2; \quad i,k=1,2,3) \label{L10}\\
{L}^{0,1}_{kli}(h) &=& 4 h  \sum_{m=1}^2 H^0_{imk} J^{1,0}_{ml}+\delta_{ki} \tau_{l}^{0,1}(h) \quad (i,l=1,2; \quad k=1,2,3) \label{L01}\\
L_{kl3}^{0,1}(h)&=&4 h D^1_{lk} +\delta_{k3}\tau^{0,1}_l(h)  \quad (k=1,2,3; \quad l=1,2)  \quad (D^1_{l3}=0)\label{L3_01}\\
{L}^{2,0}_{kli}(h) &=& {L}^0_{kli} + \dfrac{4}{h} J^{0,0}_{3l}\delta_{ki} \quad (l=1,2; \quad i,k=1,2,3) \label{L20}\\
{L}^{1,1}_{kli}(h) &=& 4 h  \sum_{m=1}^2 H^0_{imk} J^{1,0}_{lm}+\delta_{ki} \tau_{l}^{1,2}(h) \quad (i,l=1,2; \quad k=1,2,3) \label{L11}\\
L_{kl3}^{1,1}(h)&=&4 h D^1_{lk} +\delta_{k3}\tau^{1,2}_l(h)  \quad (k=1,2,3; \quad l=1,2) \quad  (D^1_{l3}=0)\label{L3_11}\\
{L}^{0,2}_{kli}(h) &=& 2  \sum_{m=1}^2  H^0_{imk}\iota^{2,1}_{lm}(h)  +\delta_{ki} \tau_{l}^{0,2}(h) \quad (i,l=1,2; \quad k=1,2,3) \label{L02}\\
{L}^{0,2}_{kl3}(h) &=& 2  \sum_{m=1}^2  \left(\vec{a}_3 \cdot \dfrac{
\partial \vec{a}_{k}  }{\partial \xi_m}\right) \iota^{2,1}_{lm}(h)  +\delta_{k3} \tau_{l}^{0,2}(h) \quad (l=1,2; \quad k=1,2,3) \label{L3_02}\\
{L}^{3,0}_{kli}(h) &=& {L}^0_{kli} + \dfrac{6}{h} J^{0,0}_{3l}\delta_{ki} \quad (l=1,2; \quad i,k=1,2,3) \label{L30}\\
{L}^{2,1}_{kli}(h) &=& 4 h  \sum_{m=1}^2 H^0_{imk} J^{1,0}_{lm}+\delta_{ki} \tau_{l}^{2,3}(h) \quad (i,l=1,2; \quad k=1,2,3) \label{L21}\\
{L}^{1,2}_{kli}(h) &=& 2  \sum_{m=1}^2  H^0_{imk}\iota^{2,1}_{lm}(h)  +\delta_{ki} \tau_{l}^{1,3}(h) \quad (i,l=1,2; \quad k=1,2,3) \label{L12}\\
{L}^{0,3}_{kli}(h) &=& 2  \sum_{m=1}^2  H^0_{imk}\iota^{3}_{lm}(h)  +\delta_{ki} \tau_{l}^{0,3}(h) \quad (i,l=1,2; \quad k=1,2,3) \label{L03}\\
\bar{S}^0_{ik}(h)&=& {S}^0_{ik} -  \dfrac{A^1}{A^0} D^0_{ik} - \dfrac{1}{h} \sum_{l=1}^2H^0_{ilk} 
J^{0,0}_{3l} \nonumber \\
&&{}-\dfrac{1}{h \sqrt{ A^0} } \left[ hI D^0_{ik}- \left( \alpha^0_i \dfrac{\partial \vec{a}_{k}}{\partial
\xi_1} + \beta^0_i \dfrac{\partial \vec{a}_{k}}{\partial
\xi_2}\right)\cdot 
\vec{\eta}(h)\right] \label{barS}\\
&=&{P}_{ik}^0+\chi(h)^0_{ik} \quad
(i=1,2; \quad k=1,2,3) \label{barSPchi}
\\
{{S}^{1,0}_{ik}}(h) &=& {S}^0_{ik}+\dfrac{2}{h}\sum_{m=1}^2 H^0_{1mk}
J^{0,0}_{3m}+
\delta_{ik} \phi^{1}(h)\quad
(i=1,2; \quad k=1,2,3)\label{S10}\\
S^{1,0}_{3k}(h) &=&{S}^0_{3k}+ \dfrac{2}{h}D^0_{3k} + \delta_{3k} \phi^{1}(h) \quad (k=1,2,3) \label{S3_10}\\
S_{ik}^{0,1}(h)&=&\sum_{l=1}^2 \left[ 2h \sum_{m=1}^2 \left( \alpha^0_i \left( \vec{a}_1 \cdot \dfrac{\partial^2 \vec{a}_k}{\partial \xi_l \partial \xi_m}\right) +  \beta^0_i \left( \vec{a}_2 \cdot \dfrac{\partial^2 \vec{a}_k}{\partial \xi_l \partial \xi_m}\right)\right)J^{1,0}_{lm}
 +    H^0_{ilk} \tau_l^{0,1}(h) \right]\nonumber\\
&& (i=1,2; \quad k=1,2,3) \label{S01}\\
S_{3k}^{0,1}(h)&=&\sum_{l=1}^2\left[ 2 h  \sum_{m=1}^2 \left(  \vec{a}_3 \cdot\dfrac{
\partial^2 \vec{a}_{k} }{\partial \xi_l \partial \xi_m}
\right) J^{1,0}_{lm}+ \left(\vec{a}_3 \cdot \dfrac{
\partial \vec{a}_{k}  }{\partial \xi_l}\right) \tau^{0,1}_l(h)  \right] \quad (k=1,2,3)\label{S3_01}\\
{{S}^{2,0}_{ik}}(h) &=& {S}^0_{ik}+\dfrac{4}{h}\sum_{m=1}^2 H^0_{imk}
J^{0,0}_{3m}+
\delta_{ik} \phi^{2,2}(h)\quad
(i=1,2; \quad k=1,2,3)\label{S20}\\
S^{2,0}_{3k}(h) &=&{S}^0_{3k}+ \dfrac{4}{h}D^0_{3k} + \delta_{3k} \phi^{2,2}(h) \quad (k=1,2,3) \label{S3_20}\\
S_{ik}^{1,1}(h)&=&\sum_{l=1}^2 \left[ 2h \sum_{m=1}^2 \left( \alpha^0_i \left( \vec{a}_1 \cdot \dfrac{\partial^2 \vec{a}_k}{\partial \xi_l \partial \xi_m}\right) +  \beta^0_i \left( \vec{a}_2 \cdot \dfrac{\partial^2 \vec{a}_k}{\partial \xi_l \partial \xi_m}\right)\right)J^{1,0}_{lm}
 +    H^0_{ilk} \tau_l^{1,2}(h) \right]\nonumber\\
&+& \delta_{ik}\phi^{1,2}(h) \quad (i=1,2; \quad k=1,2,3) \label{S11}\\
S_{3k}^{1,1}(h)&=&\sum_{l=1}^2\left[ 2 h  \sum_{m=1}^2 \left(  \vec{a}_3 \cdot\dfrac{
\partial^2 \vec{a}_{k} }{\partial \xi_l \partial \xi_m}
\right) J^{1,0}_{lm}+ \left(\vec{a}_3 \cdot \dfrac{
\partial \vec{a}_{k}  }{\partial \xi_l}\right) \tau^{1,2}_l(h)  \right]
+ \delta_{k3} \phi^{1,2}(h) \nonumber\\
&& (k=1,2,3) \label{S3_11}\\
S_{ik}^{0,2}(h)&=&\sum_{l=1}^2 \left[ \sum_{m=1}^2 \left( \alpha^0_i \left( \vec{a}_1 \cdot \dfrac{\partial^2 \vec{a}_k}{\partial \xi_l \partial \xi_m}\right) +  \beta^0_i \left( \vec{a}_2 \cdot \dfrac{\partial^2 \vec{a}_k}{\partial \xi_l \partial \xi_m}\right)\right) \iota^{2,1}_{lm}(h)
 +    H^0_{ilk} \tau_l^{0,2}(h) \right]\nonumber\\
&&  (i=1,2; \quad k=1,2,3) \label{S02}\\
S_{3k}^{0,2}(h)&=&\sum_{l=1}^2 \left[ \sum_{m=1}^2\left( \vec{a}_3 \cdot \dfrac{\partial^2 \vec{a}_k}{\partial \xi_l \partial \xi_m}\right) \iota^{2,1}_{lm}(h)
 + \left(\vec{a}_3 \cdot \dfrac{
\partial \vec{a}_{k}  }{\partial \xi_l}\right)   \tau_l^{0,2}(h) \right]\nonumber\\
&& (k=1,2,3) \label{S3_02}\\
{{S}^{3,0}_{ik}}(h) &=& {S}^0_{ik}+\dfrac{6}{h}\sum_{m=1}^2 H^0_{imk}
J^{0,0}_{3m}+
\delta_{ik} \phi^{3,3}(h)\quad
(i=1,2; \quad k=1,2,3)\label{S30}\\
S_{ik}^{2,1}(h)&=&\sum_{l=1}^2 \left[ 2h \sum_{m=1}^2 \left( \alpha^0_i \left( \vec{a}_1 \cdot \dfrac{\partial^2 \vec{a}_k}{\partial \xi_l \partial \xi_m}\right) +  \beta^0_i \left( \vec{a}_2 \cdot \dfrac{\partial^2 \vec{a}_k}{\partial \xi_l \partial \xi_m}\right)\right)J^{1,0}_{lm}
 +    H^0_{ilk} \tau_l^{2,3}(h) \right]\nonumber\\
&+& \delta_{ik}\phi^{2,3}(h) \quad (i=1,2; \quad k=1,2,3) \label{S21}\\
S_{ik}^{1,2}(h)&=&\sum_{l=1}^2 \left[ \sum_{m=1}^2 \left( \alpha^0_i \left( \vec{a}_1 \cdot \dfrac{\partial^2 \vec{a}_k}{\partial \xi_l \partial \xi_m}\right) +  \beta^0_i \left( \vec{a}_2 \cdot \dfrac{\partial^2 \vec{a}_k}{\partial \xi_l \partial \xi_m}\right)\right) \iota^{2,1}_{lm}(h)
 +    H^0_{ilk} \tau_l^{1,3}(h) \right]\nonumber\\
&+& \delta_{ik}\phi^{1,3}(h) \quad  (i=1,2; \quad k=1,2,3) \label{S12}\\
S_{ik}^{0,3}(h)&=&\sum_{l=1}^2 \left[ \sum_{m=1}^2 \left( \alpha^0_i \left( \vec{a}_1 \cdot \dfrac{\partial^2 \vec{a}_k}{\partial \xi_l \partial \xi_m}\right) +  \beta^0_i \left( \vec{a}_2 \cdot \dfrac{\partial^2 \vec{a}_k}{\partial \xi_l \partial \xi_m}\right)\right) \iota^{3}_{lm}(h)
 +    H^0_{ilk} \tau_l^{0,3}(h) \right]\nonumber\\
&&  (i=1,2; \quad k=1,2,3) \label{S03}\\
\vec{\eta}(h)&=&
\left[\dfrac{\partial h}{\partial \xi_2} (\vec{a}_1 \times
\vec{a}_3) + h \left(\vec{a}_1 \times \dfrac{\partial
\vec{a}_3}{\partial \xi_2} + \dfrac{\partial
\vec{a}_3}{\partial \xi_1} \times
\vec{a}_2\right) + \dfrac{\partial h}{\partial
\xi_1} (\vec{a}_3\times \vec{a}_2)\right]\label{eta}\\
\iota_{lm}^{2,1}(h)&=& h^2 (2 J^{2,0}_{lm}+ J^{1,1}_{lm}) \quad (l=1,2; \quad m=1,2) \label{iota21}\\
\iota^3_{lm}(h)&=&h^3(2J^{3,0}_{lm} +J^{2,1}_{lm}+ J^{2,1}_{ml} ) \quad (l=1,2; \quad m=1,2) \label{iota3}\\
\kappa(h)^0_i&=&  \sum_{l=1}^2 \dfrac{
\partial }{ \partial \xi_l}\left( \dfrac{\partial \vec{X}}{\partial t} \cdot \vec{a}_3\right) \left(L^0_{3li}-  \dfrac{A^1}{A^0}
 J^{0,0}_{il}\right)
-\dfrac{J^{0,0}_{3i}}{h^2}
   \dfrac{\partial h}{\partial t}  - \dfrac{3}{h}\sum_{k=1}^2 J^{0,0}_{ki} \dfrac{\partial^2 h}{\partial
\xi_k \partial t}\nonumber\\
&+&{}
\left( \dfrac{\partial X}{\partial t} \cdot \vec{a}_3\right) \left[{S}^0_{i3} -\dfrac{1}{h}  \sum_{l=1}^2 J^{0,0}_{3l}
H^0_{il3} 
 -\dfrac{1}{h} \dfrac{J^{0,0}_{li}}{\sqrt{A^0} } \left(  \dfrac{\partial \vec{a}_{3}}{\partial
\xi_l}\cdot 
\vec{\eta}(h)\right)  \right] \quad (i=1,2)\label{kappa} \\
{\tau}_{l}^{0,1}(h)&=&\sum_{m=1}^2 \dfrac{\partial h}{\partial \xi_m}J^{1,0}_{lm} + h K^{1,0}_{l} + h K^{0,1}_{l}+J^{1,0}_{3l}\quad (l=1,2)\label{tau1}\\
{\tau}_{l}^{1,2}(h)&=&\sum_{m=1}^2\dfrac{\partial h}{\partial \xi_m}J^{1,0}_{lm} + h K^{1,0}_{l} + h K^{0,1}_{l}+ 5 J^{0,1}_{3l}\quad (l=1,2)\label{tau12}\\
\tau_{l}^{0,2}(h)&=&h \left[\sum_{m=1}^2   \dfrac{\partial h}{\partial \xi_m}\left( 2 J^{2,0}_{lm}+  J^{1,1}_{lm} \right)+  h \left(K^{2,0}_{l} +K^{1,1}_{l}+K^{0,2}_{l}\right) + J^{1,1}_{3l} +2 J^{2,0}_{l3}\right]\nonumber \\
&&(l=1,2) \label{tau02}
\\
\tau^{2,3}_l(h)&=&\sum_{m=1}^2 \dfrac{\partial h}{\partial \xi_m} J^{1,0}_{lm} + h \left(K^{1,0}_{l} + K^{0,1}_{l}\right) + 
4J^{1,0}_{3l} + 5 J^{0,1}_{3l} \quad (l=1,2)   \label{tau23}\\
 \tau^{1,3}_l(h)&=&  h \left[\sum_{m=1}^2  \dfrac{\partial h}{\partial \xi_m}(2 J^{2,0}_{lm} +   J^{1,1}_{lm})+  h \left( K^{2,0}_{l}  +   K^{1,1}_{l} +  K^{0,2}_{l}\right) + 3  J^{1,1}_{3l} +  4    
J^{0,2}_{3l}+ 2 
  J^{2,0}_{3l}\right]\nonumber\\
&&(l=1,2) \label{tau13}\\
 \tau^{0,3}_l(h)&=& h^2\left[\sum_{m=1}^2 \dfrac{\partial h}{\partial \xi_m} (3 J^{3,0}_{lm} +  2  J^{2,1}_{lm} + J^{1,2}_{lm}) + h \left(  K^{3,0}_{l}+  K^{2,1}_{l}+  K^{1,2}_{l}+   K^{0,3}_{l}  \right)\right.\nonumber\\
& +& \left.3
 J^{0,3}_{3l}
+2  J^{1,2}_{3l} +
J^{2,1}_{3l}\right] \quad (l=1,2) \label{tau03}\\
\phi^{1}(h)&=& \dfrac{1}{h}K^{0,0}_{3}+\sum_{m=1}^2\left( H^1_{mm3}  -\dfrac{2}{h^2}\dfrac{\partial h}{\partial \xi_m} J^{0,0}_{3m} \right)\label{phi}\\
\phi^{1,2}(h)&=& \sum_{m=1}^2 \left( h H^2_{mm3} -\dfrac{1}{h}  \dfrac{\partial h }{\partial \xi_m}   J^{1,0}_{m3}\right) +K^{1,0}_{3} +K^{0,1}_{3}+ \dfrac{3}{h} J^{1,0}_{33}  \label{phi12}\\
\phi^{2,2}(h)&=& 2\sum_{m=1}^2 \left( H^1_{mm3} - \dfrac{1}{h^2} \dfrac{\partial h}{\partial \xi_m} J^{0,0}_{3m}  \right){+ \dfrac{4}{h^2}J^{0,0}_{33}}+ \dfrac{2}{h}K^{0,0}_{3}\label{phi22}\\
\phi^{3,3}(h)&=& 3\left[ \sum_{m=1}^2  \left(H^1_{mm3}  - \dfrac{1}{h^2} \dfrac{\partial h}{\partial \xi_m} J^{0,0}_{3m} \right) + \dfrac{1}{h} K^{0,0}_{3 }+ \dfrac{3}{h^2}  J^{0,0}_{33}\right] \label{phi33}\\
\phi^{2,3}(h)&=& 2 \left[\sum_{m=1}^2 \left( h  H^2_{mm3} -\dfrac{1}{h} \dfrac{\partial h}{\partial \xi_m} J^{0,1}_{3m}  \right)+  K^{0,1}_{3} +K^{1,0}_{3} +\dfrac{5}{h} J^{1,0}_{33}\right]\label{phi23}\\
\phi^{1,3}(h)&=& \sum_{m=1}^2 \left( \dfrac{\partial h }{\partial \xi_m}
(J_{3m}^{2,0} - J^{0,2}_{3m} ) +  h^2 H^3_{mm3}\right) + h (K^{2,0}_{3}+   K^{1,1}_{3} +  K^{0,2}_{3}) \nonumber\\
&&{} + 4 J^{2,0}_{33} + 2 J^{1,1}_{33} \label{phi13}\\
\chi(h)_{ik}^0&=&\dfrac{1}{h}   \left\{  \dfrac{\partial h}{\partial
\xi_1} \left[\sum_{l=1}^2H^0_{ilk}\alpha_l^0 -\dfrac{1}{\sqrt{A^0}}(\vec{a}_3\times \vec{a}_2)\cdot \left( \alpha_i^0
\dfrac{
\partial \vec{a}_{k}}{\partial \xi_1}+\beta_i^0 \dfrac{
\partial \vec{a}_{k}}{\partial \xi_2}\right)\right]\right.\nonumber\\
&+&\left. \dfrac{\partial h}{\partial \xi_2} \left[ \sum_{l=1}^ 2H^0_{ilk}\beta_l^0 -\dfrac{1}{\sqrt{A^0}} (\vec{a}_1 \times
\vec{a}_3) \cdot \left( \alpha_i^0
\dfrac{
\partial \vec{a}_{k}}{\partial \xi_1}+ \beta_i^0
\dfrac{
\partial \vec{a}_{k}}{\partial \xi_2} \right)  \right]\right\} \nonumber\\
&& (i=1,2; \quad k=1,2,3)\label{chiik}\\
\psi (h)_{ijl}^{0}&=& \dfrac{1}{h} \left [ \left ( \alpha _{l}^{0}
\dfrac{\partial h}{\partial \xi _{1}} +\beta _{l}^{0}
\dfrac{\partial h}{\partial \xi _{2}} \right ) \delta _{ij} +
\dfrac{\partial h}{\partial \xi _{j}}J^{0,0}_{li}\right ] \quad (i,j,l=1,2)
\label{psi}
\end{eqnarray}
where $\delta_{ij}$ is the Kronecker Delta.

\end{document}